\documentclass{article}
\usepackage[english]{babel}

\usepackage{delarray}
\usepackage{amsmath, latexsym, amsfonts, amssymb, amsthm, amscd}
\usepackage{graphicx}
\newtheorem {lemme} {{\bf Lemma}} [section]
\newtheorem {theoreme} {{\bf Theorem}} [section]
\newtheorem {proposition} {{\bf Proposition}} [section]

\newtheorem {remarque} {{\bf Remark}} [section]

\newcommand{\tr}{\operatorname{tr}}

\newcommand{\bX}{{\bf X}}
\newcommand{\bW}{{\bf W}}
\newcommand{\bA}{{\bf A}}
\newcommand{\bB}{{\bf B}}

\newcommand{\bM}{{\bf M}}

\newcommand{\bQ}{{\bf Q}}
\newcommand{\Tr}{\operatorname{Tr}}
\newcommand{\1}{1\!\!{\sf I}_{\Omega_N(\Lambda_N)}}
\newcommand{\deux}{1\!\!{\sf I}_{\Omega_N^{(2)}}}

\newcommand{\E}{{\mathbb E}}

\newcommand{\PP}{{\mathbb P}}
\newcommand{\R}{{\mathbb R}}
\newcommand{\C}{{\mathbb C}}

\newcommand{\vers}{\mathop{\longrightarrow}}
\newcommand{\lmax}{\lambda_1(\bM_N)}

\numberwithin{equation}{section}
\title{ Central limit theorems for eigenvalues  of deformations of Wigner matrices\footnote{This work was partially supported by the {\emph Agence Nationale de la
Recherche} grant ANR-08-BLAN-0311-03.}
}
\author{M. Capitaine\thanks{CNRS, Institut de Math\'ematiques de Toulouse, Equipe de Statistique et Probabilit\'es,  F-31062 Toulouse Cedex 09. E-mail: mireille.capitaine@math.univ-toulouse.fr } ,
 \hspace{.1cm}C. Donati-Martin\thanks{UPMC Univ Paris 06 and CNRS, UMR 7599, Laboratoire de Probabilit\'es et Mod\`eles Al\'eatoires,  Case 188, 4 palce Jussieu, F-75252 Paris Cedex 05. E-mail: catherine.donati@upmc.fr }
  \hspace{.1cm}   and 
D. F\'eral\thanks{Institut de Math\'ematiques de Bordeaux, Universit\'e  Bordeaux 1, 351 Cours de la Lib\'eration, F-33405 Talence Cedex. E-mail: delphine.feral@math.u-bordeaux1.fr} }

\date{}

\begin{document}
\maketitle
\begin{abstract} In this paper, we  study  the fluctuations of the extreme eigenvalues  of a spiked finite rank deformation of a Hermitian (resp. symmetric) Wigner matrix when  these eigenvalues separate from the bulk.
We exhibit quite general situations that will give rise to universality or non universality of the fluctuations, according to the delocalization or localization of the eigenvectors of the perturbation. Dealing with the particular case of a spike with multiplicity one, 
we also establish  a necessary and sufficient condition on the associated normalized eigenvector  so that the fluctuations of the corresponding eigenvalue of the deformed model are universal.
\end{abstract}
\begin{center}
R\'esum\'e
\end{center}
Dans ce papier, nous \'etudions les fluctuations des valeurs propres extr\'emales d'une matrice de Wigner hermitienne (resp. sym\'etrique) d\'eform\'ee par une perturbation de rang fini dont les  valeurs propres  non nulles sont fix\'ees, dans le cas o\`u ces  valeurs propres extr\'emales  se d\'etachent du reste du spectre.
Nous d\'ecrivons des situations g\'en\'erales d'universalit\'e ou de non-universalit\'e des fluctuations correspondant au caract\`ere localis\'e ou d\'elocalis\'e des vecteurs propres de la perturbation. Lorsque l'une des valeurs propres de la perturbation est de multiplicit\'e un,  nous \'etablissons de plus une condition n\'ecessaire et suffisante sur le vecteur propre associ\'e pour que les fluctuations de la valeur propre correspondante du mod\`ele d\'eform\'e soient universelles.

\vspace{.3cm}
\noindent {\it Mathematics Subject Classification (2010)}: 60B20, 15A18, 60F05 . \\
\noindent {\it Key words}:  Random matrices, deformed Wigner matrices,  extremal eigenvalues,
fluctuations, localized eigenvectors, universality.

\section{Introduction}
Adding a finite rank perturbation to a GUE matrix, S. P\'ech\'e \cite{Pe} pointed out a sharp phase transition phenomenon :  according to the largest eigenvalue of the  perturbation, the largest eigenvalue of the perturbed matrix   should either  stick to the bulk and fluctuate
according to the Tracy-Widom (or generalized Tracy-Widom ) law or should be extracted away from the bulk and have then  fluctuations of Gaussian nature. In the lineage of this work, in a previous paper
 \cite{CDF}, we have studied the limiting behavior of extremal eigenvalues of finite rank deformations of Wigner matrices. We established their almost sure convergence. The limiting values depend only on the spectrum of the deformation $A_N$ and on the variance of the distribution of the entries of the Wigner matrix. 
 On the contrary the fluctuations of these eigenvalues strongly depend on the eigenvectors of $A_N$. 
Indeed, in the particular case of a rank one diagonal deformation whose non-null eigenvalue is large enough, we established a central limit theorem for the largest eigenvalue which deviates from the rest of the spectrum
and  proved that the fluctuations of the largest eigenvalue vary with the particular distribution of the entries of the Wigner matrix. Thus,  this fluctuations result differs from that of the full rank one deformation case investigated in \cite{FK} and \cite{FePe} since this latter case exhibited universal limiting distributions. \\
Let us recall these results in the complex setting, having in mind that similar results hold in the real symmetric case. In the following, given an arbitrary Hermitian matrix $M$ of size $N$, we will denote by $\lambda_1(M) \geq \cdots \geq \lambda_N(M)$
its $N$ ordered eigenvalues; we will denote the centered gaussian distribution with variance $v$ by ${\cal N}(0,v)$. \\
The random matrices under study are complex Hermitian matrices $(\bM_N)_N$  defined on a probability space $(\Omega, {\cal F}, \PP)$ such that
\begin{equation}
\bM_N=\frac{\bW_N}{\sqrt{N}}+\bA_N.
\end{equation}
 $\bA_N$ is a $N \times N$ deterministic Hermitian 
  matrix of fixed finite rank and whose spectrum does not depend on $N$. The matrix $\bW_N$ is a $N \times N$ Wigner Hermitian  matrix such that the $N^2$ random variables
$(W_N)_{ii}$, $\sqrt{2} \Re e((W_N)_{ij})_{i<j}$, $\sqrt{2}
\Im m ((W_N)_{ij})_{i<j}$ 
are independent identically distributed with
a centered distribution $\mu$ of variance $\sigma^2$.

\noindent As the rank of the $\bA_N$'s is assumed to be finite, the Wigner
Theorem is still satisfied for the Deformed Wigner model $(\bM_N)_N$ (cf. Lemma 2.2 of \cite{Bai}): the spectral measure $\frac{1}{N} \sum_{i=1}^N
\delta_{\lambda_i(\bM_N)}$ of $\bM_N$ converges a.s.
towards the semicircle law $\mu _{sc}$ whose density is given by
\begin{equation}\label{scl}
\frac{d \mu_{sc}}{dx}(x)=  \frac{1}{2 \pi \sigma^2} \sqrt{4\sigma^2
- x^2} \, 1 \hspace{-.20cm}1_{[-2\sigma,2 \sigma]}(x).
\end{equation}

\noindent  When $\bA_N \equiv 0$, it is well-known that once $\mu$ has a finite fourth moment, the first largest (resp. last smallest) eigenvalues of the rescaled Wigner matrix $\bW_N / \sqrt N$ tend almost surely to the right (resp. left)-endpoint $2 \sigma$ (resp. $-2 \sigma$) of the semicircle support (cf.
\cite{Bai}). The corresponding fluctuations, which have been first obtained by Tracy and Widom \cite{TW}     in the Gaussian case and then extended by Soshnikov \cite{So} for any symmetric probability measure $\mu$ having subgaussian moments, are governed by the so-called Tracy-Widom distributions. Note that the exponential decay condition (with symmetry assumption)  has been replaced by a finite number of moments in  \cite{Ru}, \cite{Kho}.
Under the subexponential decay assumption, the symmetry assumption on $\mu$ in \cite{So} was replaced in \cite{TV} by the vanishing third moment condition and very recently, Erd\"{o}s, Yau and Yin \cite{EYY} proved the edge universality under the subexponential decay assumption alone. \\
Let us describe how the asymptotic behavior of the extremal
eigenvalues of the perturbed Wigner matrix may be affected by the perturbation by considering the particular case of a rank one perturbation $\bA_N$ with non-null eigenvalue $\theta$. For a large class of  probability measures $\mu$, it turns out that
the largest eigenvalue $\lambda_1(\bM_N)$ still tends to the right-endpoint $2 \sigma$ if $\theta \leq \sigma$ whereas
 $\lambda_1(\bM_N)$ jumps above the bulk to $\rho_\theta = \theta + \frac{\sigma^2}{\theta}$ if $\theta > \sigma$.
This was proved by P\'ech\'e in her pionnering work \cite{Pe} when $\mu$ is gaussian, extended in \cite{FePe} when $\mu$ is symmetric and has subgaussian moments but in the particular case of the full rank one deformation $\bA_N$ given by \begin{equation}\label{fullrank}(\bA_N)_{ij} = {\theta}/{N} \mbox{ ~ for all~} 1\leq i,j \leq N \end{equation}  and
finally established in \cite{CDF} for general  $\bA_N$ when $\mu$ is symmetric and  satisfies a Poincar\'e inequality. 

\noindent Moreover, considering the perturbation matrix defined by (\ref{fullrank}),
 F\'eral and P\'ech\'e \cite{FePe} proved that the fluctuations of $\lambda_1(\bM_N)$ are the same as in the gaussian setting investigated in \cite{Pe} and in this sense are universal. Here is their result when $\theta> \sigma$:
\begin{proposition}\label{TCLFePe}
If $\mu$ is symmetric and has subgaussian moments
$$\sqrt{N} (\lambda_1(\bM_N) - \rho_\theta) \overset{\mathcal
  L}{\longrightarrow}  {\cal
  N}(0, \sigma^2_{\theta})$$
  where $\sigma_\theta = \sigma \sqrt{ 1-\frac{\sigma^2}{\theta^2}}$.
  \end{proposition}
\noindent   The proof of this result relies on the computations of moments of $\bM_N$ of high order (depending on $N$) and the knowledge of the fluctuations in the Gaussian case, established by P\'ech\'e \cite{Pe}. 
  
 \noindent  On the other hand, for the strongly localized perturbation matrix of rank 1 given by 
  $$\bA_N={\rm{diag}}(\theta, 0, \cdots,0)$$ \noindent with  $\theta > \sigma$, we  proved in \cite{CDF} that the fluctuations of $\lambda_1(\bM_N)$ vary with the particular distribution of the entries of the Wigner matrix so that this phenomenon can be seen as an example of a non universal behavior :
\begin{proposition}  {\label{ThmFluctuations}}
Let $\mu$ be symmetric and satisfy a Poincar\'e inequality.
Define
\begin{eqnarray}\label{defvtheta} 
c_{\theta}= \frac{\theta^2}{\theta^2- \sigma^2} \quad \text{and} \quad  v_{\theta}= \frac{1}{2}\Big ( \frac{m_4 - 3 \sigma^4}{\theta^2} \Big )  +
\frac{\sigma ^4}{\theta^2-\sigma^2}
\end{eqnarray}
where  $m_4:= \int x^4 d\mu(x)$. Then
\begin{eqnarray}\label{FluctDef}
 c_{\theta} \sqrt{N} \Big (\lmax -\rho_{\theta} \Big ) \overset{\mathcal
  L}{\longrightarrow}  \Big{\{} \mu \ast {\cal
  N}(0, v_{\theta})\Big{\}}.
\end{eqnarray}
\end{proposition}

In the present paper, we consider perturbations $\bA_N$ of higher rank of Wigner matrices associated to some symmetric probability measure $\mu$
satisfying a Poincar\'e inequality. The a.s. convergence of the extreme eigenvalues has already been described in \cite{CDF} (see 
Theorem \ref{defrhotheta} below). Whenever the largest eigenvalues of  $\bM_N$ are extracted away from the bulk, we describe their fluctuations which depend on the localization of the eigenvectors of $\bA_N$, as already seen in the above rank 1 examples. We investigate two quite general situations for which we exhibit a phenomenon of different nature. To explain this, let us focus on the largest eigenvalue $\theta_1$ of $\bA_N$. We assume that $\theta_1 > \sigma$ so that the  largest eigenvalues of  $\bM_N$ converges a.s towards $\rho_{\theta_1} = \theta_1 + \frac{\sigma^2}{\theta_1}> 2\sigma$ 
. \\
First, when the  eigenvectors associated to the largest eigenvalue $\theta_1$ of $\bA_N$ 
 are localized,
 we establish that the
 limiting distribution in the fluctuations of
 $  \lambda_{i}(\bM_N),$ $ 1 \leq i \leq k_1$,
around $ \rho _{\theta_1} $
is not universal and we give it explicitely in terms of these eigenvectors and of the distribution of the entries of the Wigner matrix, see Theorem \ref{CLT}. \\ 
Secondly, if the eigenvectors are sufficiently delocalized, we establish the universality of the fluctuations of $ \lambda_{i}(\bM_N),$ $ 1 \leq i \leq k_1$, see Theorem \ref{CLTkN} .  \\ 
Actually, in the rank one case,  this study allows us to exibit a necessary and sufficient condition on a normalized eigenvector of $A_N$ associated to the largest eigenvalue $\theta_1$ for the universality of the fluctuations (see Theorem \ref{rank1} below). 
Moreover if such an   eigenvector of $\bA_N$  is not localized  but does not   satisfy  the  criteria of  universality,
 the largest eigenvalue of $\bM_N$ may  fluctuate  according to a mixture of $\mu$ and normal distributions generalizing (\ref{FluctDef}). We will describe some of such intermediate situations.\\
We detail the definition of localization/delocalization and these results in the Section 3.

\noindent Our approach is close to the proof of (\ref{FluctDef}), with more involved computations  and are in the spirit of the works of  \cite{Pa}
and  \cite{BBPbis}. It is valid in both real and complex settings.
Actually, we assume that the eigenvectors associated to the largest eigenvalues  of $\bA_N$ belong to a subspace generated by $k (=k(N))$ canonical vectors of $\C^N$ and the method requires that $N-k \vers \infty$ (and even $\frac{k}{\sqrt{N}} \vers 0$). In particular, this approach does not cover the case of Proposition \ref{TCLFePe}  studied by \cite{FePe} where $k=N$. \\

The Deformed Wigner matrix model may be seen as the additive analogue of the
spiked population models. These are random sample covariance matrices
$(S_N)_N$ defined by $S_N=\frac{1}{N} Y_N^* Y_N$ where $Y_N$ is a $p \times N$ complex (resp. real) matrix (with $N$ and $p=p_N$ of the same order
as $N \to \infty$) whose entries satisfy first four moments conditions; the sample column
vectors are assumed to be i.i.d, centered and of covariance matrix a
deterministic Hermitian (resp. symmetric) matrix ${\Sigma}_p$ having all but finitely
many eigenvalues equal to one. In their pioneering article on that topic \cite{BBP}, Baik-Ben Arous-P\'ech\'e  pointed out a phase  transition phenomenon for the fluctuations of the  largest eigenvalue of $S_N$ according to the largest eigenvalue of ${\Sigma}_p$, in the complex Gaussian setting; their results were extended in \cite{Pa} to the real case when the largest eigenvalue of ${\Sigma}_p$ is simple and sufficiently larger than 1 and in \cite{O} to singular Wishart matrices.
In the non Gaussian case, the fluctuations of the extreme eigenvalues 
have been recently studied by   Bai-Yao  \cite{BY2} and F\'eral-P\'ech\'e \cite{FePe2}.\\

The paper is organized as follows. In Section 2, we present the matricial models under study and the notations that will be used throughout the paper. In Section 3,
 we present the main results of this paper. We  give a summary of our approach in Section 4. 
Section 5 is
devoted to the proof of Theorem \ref{CLT}, Theorem \ref{CLTkN} and Theorem \ref{rank1}.  
 Finally, we recall some basic facts on matrices, a CLT for random sesquilinear forms and prove some technical results
in an Appendix. \\

\section{Model and notations}
 The random matrices under study are complex Hermitian (or real symmetric) matrices $(\bM_N)_N$  defined on a probability space $(\Omega, {\cal F}, \PP)$ such that
\begin{equation}{\label{defM}}
\bM_N=\frac{\bW_N}{\sqrt{N}}+\bA_N
\end{equation}
where the matrices $\bW_N$ and $\bA_N$ are defined as follows:
\begin{itemize}
\item[{\bf{(i)}}] $\bW_N$ is a $N \times N$ Wigner Hermitian (resp. symmetric) matrix such that the $N^2$ random variables
$(W_N)_{ii}$, $\sqrt{2} \Re e((W_N)_{ij})_{i<j}$, $\sqrt{2}
\Im m ((W_N)_{ij})_{i<j}$ (resp. the $\frac{N(N+1)}{2}$ random variables $\frac{1}{\sqrt{2}} (W_N)_{ii}$,
$(W_N)_{ij}$, $i<j$) are independent identically distributed with
a symmetric distribution $\mu$ of variance $\sigma^2$ and satisfying a
Poincar\'e inequality; the latter condition means that
 there exists a positive constant $C$ such that for any
${\cal C}^1$ function $f: \R \rightarrow \C$  such that $f$ and
$f'$ are in $L^2(\mu)$,
$$\mathbf{V}(f)\leq C\int \vert f' \vert^2 d \mu ,$$
\noindent with $\mathbf{V}(f) = \mathbb{E}(\vert
f-\mathbb{E}(f)\vert^2)$. \\
Note that when $\mu$ is Gaussian, $W_N$ is a GU(O)E$(N \times N, \sigma^2)$ matrix.
\item[{\bf{(ii)}}] $\bA_N$ is a deterministic Hermitian (resp. symmetric) matrix of fixed finite
rank $r$ and built from a family of $J$ fixed real numbers
$\theta_1>\cdots > \theta_J$ independent of $N$ with some $j_0$ such
that $\theta _{j_0}=0$. We assume that the non-null eigenvalues
$\theta_j$ of $\bA_N$ are of fixed multiplicity $k_{j}$ (with $\sum_{j
\not= j_0} k_j=r$). 
Let $J_{+ \sigma}$  be the number of j's such that
$\theta_j > \sigma$.
We denote by  $k_{+ \sigma}:=k_1 + \cdots + k_{J_{+ \sigma}}$.
 We introduce   $k \geq k_{+ \sigma}$  as the minimal number of canonical vectors among the canonical basis $(e_i; \, i= 1, \ldots, N)$ of $\C^N$ needed to express all the eigenvectors  associated to the largest eigenvalues $\theta_1, \ldots, \theta_{J_{+ \sigma}}$  of $\bA_N$.
 Without loss of generality (using the invariance of the distribution of the Wigner matrix $\bW_N$ by conjugation by a permutation matrix), we can assume that these $k_{+ \sigma}$ eigenvectors belong to ${\rm{Vect}} (e_1, \ldots, e_k)$.  \\
 All along the paper we assume that $k\ll \sqrt{N}$.
\end{itemize}
Let us now fix $j$ such that $1 \leq j \leq J_{+ \sigma}$ and let $U_k$ be a unitary matrix  of size $k$ such that  \begin{equation} \label{hypA_N}
{\rm diag}(U_k^*, I_{N-k}) \bA_N {\rm diag}(U_k, I_{N-k})  = {\rm diag}(
\theta_j I_{k_j},   (\theta_{l}I_{k_{l}})_{ l \leq J_{+ \sigma}, l \not= j} , Z_{N-k_{+ \sigma}})
\end{equation}
where $Z_{N-k_{+ \sigma}}$ is an Hermitian matrix with eigenvalues strictly smaller than $\theta_{J_{+ \sigma}}$.\\

\noindent Define $K_j=K_j(N)$ as the minimal number of canonical vectors among $ (e_1, \ldots, e_k)$ needed to express all 
the orthonormal eigenvectors $v^{j}_i$, $ 1 \leq i \leq k_j$, of $\bA_N$ associated to $\theta_j$.  Without loss of generality, we can assume that the $v^{j}_i$ , $ 1 \leq i \leq k_j$, belong to ${\rm{Vect}} (e_1, \ldots, e_{K_j})$.
Considering now the vectors $v^i_j$ as vectors in $\mathbb C^{K_j}$, we define the ${K_j \times k_j}$ matrix
\begin{equation}\label{defUj} U_{K_j \times k_j}:= \left( v^j_1,\ldots, v^j_{k_j} \right)\end{equation} \noindent namely $U_{K_j \times k_j}$ is the upper left corner of $U_k$ of size $K_j \times k_j$. It satisfies
\begin{equation}\label{propUj} U_{K_j \times k_j} ^* U_{K_j \times k_j}=I_{k_j}. \end{equation}
All along the paper, the parameter $t$ is such that $t=4$ (resp. $t=2$) in the real (resp. complex) setting and we let $m_4:= \int x^4 d\mu(x)$. \\

\noindent Given an arbitrary Hermitian or symmetric matrix $M$ of size $N$, we will denote by $\lambda_1(M) \geq \cdots \geq \lambda_N(M)$
its $N$ ordered eigenvalues.\\

\section{Main results}
We first recall the a.s. convergence of the extreme eigenvalues.
Define
\begin{eqnarray}{\label{defrhotheta}}
\rho
_{\theta_j}=\theta_j + \frac{\sigma^2}{\theta_j}.
\end{eqnarray}
Observe that $\rho
_{\theta_j} >2 \sigma$ (resp. $<-2 \sigma$) when $\theta_j>  \sigma$
(resp. $<-\sigma$) (and $\rho
_{\theta_j} =\pm 2 \sigma$ if $\theta _j= \pm \sigma$). \\
For definiteness, we set $k_1+ \cdots +k_{j-1}:=0$ if $j=1$.
In \cite{CDF}, we have established the following universal convergence result.
\begin{theoreme}{\label{ThmASCV}} ({\bf a.s. behaviour})
Let $J_{+ \sigma}$ (resp. $J_{- \sigma}$) be the number of j's such that
$\theta_j > \sigma$ (resp. $\theta_j < -\sigma$).
\begin{itemize}
\item[{(1)}] $\quad \forall 1 \leq j \leq J_{+ \sigma}, \, \forall 1 \leq i \leq k_j, \quad  \lambda_{k_1+ \cdots +k_{j-1} +i}(\bM_N)
\longrightarrow \rho _{\theta_j}  \quad{a.s.}$
\item[\text{(2)}] $\quad \lambda_{k_1+ \cdots +k_{J_{+ \sigma}} +1}(\bM_N)
\longrightarrow 2\sigma \quad{a.s.}$
\item[\text{(3)}] $\quad \lambda_{k_1+ \cdots +k_{J-J_{- \sigma}}}(\bM_N)
\longrightarrow -2 \sigma \quad{a.s.}$
\item[\text{(4)}] $\quad \forall j \geq J-J_{- \sigma}+1, \, \forall 1 \leq i \leq k_j, \quad
\lambda_{k_1+ \cdots +k_{j-1} +i}(\bM_N) \longrightarrow \rho _{\theta_j} \quad{a.s}.$
\end{itemize}
\end{theoreme}

\noindent 
 
\subsection{ Fluctuations around $\rho_{\theta_j}$}
>From Theorem \ref{ThmASCV},
 for all $1 \leq i \leq k_j$, $\lambda_{k_1+ \cdots +k_{j-1} +i}(\bM_N)$ converges to
$ \rho _{\theta_j}$ a.s..
We shall describe their fluctuations   in the extreme two cases: \\

\noindent {\bf Case a) localization of the eigenvectors associated to $\theta_j$:} The sequence $K_j(N)$ is bounded, 
$$\sup_N K_j(N) = \tilde{K_j}$$ 
and the the upper left corner $U_{\tilde{K}_j\times k_j}$ of $U_k$ of size $\tilde{K}_j \times k_j$  converges towards some matrix $\tilde{U}_{\tilde{K_j}\times k_j}$ when $N$ goes to infinity;\\

\noindent {\bf Case b)  delocalization of the eigenvectors associated to $\theta_j$:} $K_j=K_j(N) \rightarrow
\infty$ when $N \rightarrow \infty$ and  $U_k$ satisfies
\begin{equation} \label{hypU_k}
 \max_{p=1}^{k_j}\max_{i=1}^{K_j} \vert (U_k)_{ip}\vert \longrightarrow 0 \quad \text{as } ~{ N \rightarrow \infty}. \mbox{~}
 \end{equation}

\noindent
The main results of our paper are  the following two theorems.  Let  $c_{\theta_j}$ be  defined by
\begin{equation} \label{defctheta}
c_{\theta_j} = \frac{\theta_j^2}{\theta_j^2- \sigma^2}. \end{equation}

In Case a) (which includes the particular setting of Proposition \ref{ThmFluctuations}), the fluctuations of the corresponding  rescaled largest eigenvalues of $\bM_N$ are not universal.

\begin{theoreme} \label{CLT} In Case a): the $k_j$-dimensional vector
$$ \left (c_{\theta_j} \sqrt{N} (\lambda_{k_1 + \ldots + k_{j-1} + i}( \bM_N)- \rho_{\theta_j}); \, i = 1, \ldots, k_j \right )$$
converges in distribution to $(\lambda_i( V_{k_j\times k_j}); i = 1, \ldots k_j)$ where $\lambda_i(V_{k_j\times k_j})$ are the ordered eigenvalues of the matrix  $V_{k_j\times k_j}$ of size $k_j$
defined in the following way. Let
$W_{\tilde{K}_j} $ be a Wigner matrix of size ${\tilde{K}_j} $ with distribution given by $\mu$ (cf {\bf (i)}) and $H_{\tilde{K}_j} $
 be a centered Hermitian Gaussian matrix of size ${\tilde{K}_j} $ independent of $W_{\tilde{K}_j} $ with independent entries $H_{pl}$, $p \leq l$ with variance
\begin{equation}
\left\{  \begin{array}{l} \displaystyle v_{pp} = E(H_{pp}^2) = \frac{t}{4}\Big ( \frac{m_4 - 3 \sigma^4}{\theta_j^2} \Big )
   + \frac{t}{2}\frac{\sigma^4 }{\theta_j^2-\sigma^2} ,  \, p= 1, \ldots, {\tilde{K}_j} , \\
\displaystyle v_{pl} = \mathbb{E}(|H_{pl}|^2) =
\frac{\sigma^4}{\theta_j^2-\sigma^2},  \, 1 \leq p<l \leq {\tilde{K}_j} .
\end{array} \right.  \end{equation}
Then, $V_{k_j\times k_j}$
is the  $k_j \times k_j$ matrix defined by
\begin{equation} V_{k_j\times k_j} =  \tilde{U}_{\tilde{K}_j\times k_j}^* (W_{\tilde{K}_j}  + H_{\tilde{K}_j} ) \tilde{U}_{\tilde{K}_j\times k_j} .\end{equation}
\end{theoreme}

Case b) exhibits universal fluctuations.

\begin{theoreme} \label{CLTkN} In  Case b): the $k_j$-dimensional vector
$$ \left (c_{\theta_j} \sqrt{N} (\lambda_{k_1 + \ldots + k_{j-1} + i}( \bM_N)- \rho_{\theta_j}); \, i = 1, \ldots, k_j \right )$$
converges in distribution to $(\lambda_i( V_{k_j\times k_j}); \, i = 1, \ldots k_j)$ where the matrix $V_{k_j\times k_j}$
is distributed as the GU(O)E($k_{j}\times k_{j}, \frac{\theta_j^2 \sigma^2}{\theta_j^2- \sigma^2}$).
\end{theoreme}
\begin{remarque} Note that since $\mu$ is symmetric, analogue results can be deduced from Theorem \ref{CLT} and Theorem \ref{CLTkN} dealing with the lowest eigenvalues of $\bM_N$ and the $\theta_j$ such that $\theta_j < - \sigma$.
\end{remarque}

\noindent
{\bf Example:}
$$A_N = {\rm diag}(A_p(\theta_1), \theta_2 I_{k_2}, 0_{N-p-k_2})$$
where $A_p(\theta_1)$ is a matrix of size $p$ defined by $A_p(\theta_1)_{ij}  = \theta_1 /p$, with $\theta_1, \theta_2 > \sigma$, $p \ll \sqrt{N}$.
Then $k = p + k_2$, $k_1 = 1$, $K_1 = p$, $K_2 = k_2$. For $j=1$, we are in Case a) if $p$ is bounded and in Case b) if $p= p(N) \rightarrow + \infty$. For $j=2$, we are in Case a).

\subsection{Further result for  a spike $\theta_j > \sigma$ of multiplicity  1}\label{intermediate}
Dealing with a spike $\theta_j > \sigma$ with multiplicity  1, it turns out that case b) is actually the unique situation where universality holds since we establish the following.
\begin{theoreme}\label{rank1}
If $k_j=1$, $\theta_j >\sigma$, 
then the fluctuations of $\lambda_{k_1+\cdots+k_{j-1}+1}(\bM_N)$ are universal, namely $$\sqrt{N} (\lambda_{k_1+\cdots+k_{j-1}+1}(\bM_N) - \rho_{\theta_j}) \overset{\mathcal
  L}{\longrightarrow}  {\cal
  N}(0,\frac{t}{2} \sigma^2_{\theta_j}) \mbox{~where~} \sigma_{\theta_j}=\sigma \sqrt{ 1-\frac{\sigma^2}{\theta_j^2}},$$
  if and only if \begin{equation}\label{cns}\max_{l \leq K_j} \vert (U_k)_{l 1}\vert \rightarrow 0 \mbox{~when $N \rightarrow \infty$}. \end{equation}
\end{theoreme}

Moreover, our approach allows us to     describe the fluctuations of $\lambda_{k_1+\cdots+k_{j-1}+1}(\bM_N)$ for some particular situations where the corresponding eigenvector of $\bA_N$  is not localized  but does not   satisfy  the  criteria of  universality $\max_{l \leq K_j} \vert (U_k)_{l 1}\vert \rightarrow 0$ (that is somehow for intermediate situations
between Case a) and Case b)). Let $m$ be a fixed integer number. Assume that for any $l =1,\ldots,m$,
$(U_k)_{l 1}$ is independent of $N$,   whereas $\max_{m < l \leq K_j} \vert (U_k)_{l 1}\vert \rightarrow 0$ when $N$ goes to infinity. We will prove at the end of Section \ref{preuves} that $c_{\theta_j}\sqrt{N} (\lambda_{k_1+\cdots+k_{j-1}+1}(\bM_N) - \rho_{\theta_j})$ converges in distribution towards
the mixture of $\mu$-distributed or gaussian random variables $\sum_{i,l=1}^m a_{il} \xi_{il}+ {\cal N}$ in the complex case, $\sum_{1\leq l \leq i \leq m} a_{il} \xi_{il}+ {\cal N}$
in the real case,
where $\xi_{il}, (i,l) \in \{1,\ldots,m\}^2,~ {\cal N}$ are independent random variables such that
\begin{itemize}
\item for any $(i,l) \in \{1,\ldots,m\}^2$, the distribution of $\xi_{il}$ is $\mu$;
\item $a_{il}=\left\{\begin{array}{lll}\sqrt{2} \Im (\overline{(U_k)_{l 1}}(U_k)_{i 1}) \mbox{~if~} i<l\\
\sqrt{t} \Re (\overline{(U_k)_{l 1}}(U_k)_{i 1}) \mbox{~if~} i>l\\ \sqrt{\frac{t}{2}}
\vert (U_k)_{l 1}\vert^2 \mbox{~if~} i=l; \end{array} \right.$
\item ${\cal N}$ is a centered  gaussian variable with variance 
$$ \frac{t}{4} \frac{\left[m_4 - 3 \sigma^4\right] \sum_{l=1}^m\vert (U_k)_{l 1}\vert^4}{\theta_j^2}   +
\frac{t}{2}\frac{\sigma ^4}{\theta_j^2-\sigma^2} + \frac{t}{2}\left[1- \left(\sum_{l=1}^m\vert (U_k)_{l 1}\vert^2\right)^2\right]\sigma^2.$$
\end{itemize}

\section{Sketch of the approach}

Before we proceed to the proof of Theorems \ref{CLT} and \ref{CLTkN}, let us give the sketch of our approach which are similar in both cases. To this aim, we define for any random variable $\lambda$,
\begin{equation} \label{defxi}
\xi_N(\lambda) = c_{\theta_j} \sqrt{N}(\lambda- \rho_{\theta_j})
\end{equation}
with $c_{\theta_j}$ given by \eqref{defctheta}. We also set $ \hat{k}_{j-1}:= k_1 + \ldots +  k_{j-1}$ with the
convention that $\hat{k}_{0}=0$.\\
The reasoning made in the setting of Proposition \ref{ThmFluctuations} (for which $k=k_{+\sigma} =1$) relies (following ideas previously developed in \cite{Pa} and \cite{BBPbis}) on
the writing of the rescaled eigenvalue $\xi_N( \lambda_{  1}( \bM_N))$ in terms of the resolvent of an underlying non-Deformed Wigner matrix.
The conclusion then
essentially follows from a CLT on random sesquilinear forms established by J. Baik and J. Silverstein in the Appendix of \cite{CDF} (which corresponds to the following Theorem \ref{Silver} in the scalar case). In the general case, to prove the convergence in distribution of the vector $ \big ( \xi_N( \lambda_{\hat{k}_{j-1} + i}( \bM_N)); i = 1, \ldots, k_j \big )$, we will extend, as \cite{BY2}, the previous approach  in the following sense. We will show that each of these rescaled eigenvalues is an eigenvalue of a $k_j \times k_j$ random matrix which may be expressed in terms of the resolvent of a $N-k \times N-k$ Deformed Wigner matrix whose eigenvalues do not jump asymptotically outside $[-2 \sigma;2\sigma]$; then, the matrix $V_{k_j \times k_j}$ will arise from a multidimensional CLT on random sesquilinear forms. Nevertheless, due to the multidimensional situation to be considered now, additional considerations are required. Let us give more details.\\
Consider an arbitrary random variable $\lambda$ which converges in probability towards $\rho_{\theta_j}$. Then, applying factorizations of type \eqref{calculdet}, we prove that $\lambda$ is an eigenvalue of $\bM_N$ iff $\xi_N(\lambda)$ is (on some event having probability going to 1 as $N \to \infty$) an eigenvalue of a $k_{j} \times k_{j}$ matrix ${\check X}_{k_{j},N}( \lambda)$ of the form
\begin{equation}\label{ecriture1adeX}
\check{X}_{k_{j},N}( \lambda)=  V_{k_j,N}+ R_{k_{j},N}(\lambda)
\end{equation}
where $V_{k_j,N}$ converges in distribution towards $V_{k_{j} \times k_{j}}$ and  the remaining term $R_{k_{j},N}(\lambda)$ turns out to be negligible.
Now, when $k_j>1$, since the matrix $
\check{X}_{k_{j},N}( \lambda)$ (in \eqref{ecriture1adeX}) depends on $\lambda$,
the previous reasoning with $\lambda=\lambda_{\hat k _{j-1} + i}({\bM}_N)$ for any $1 \leq i \leq k_j$ does not allow us to readily deduce  that the $k_j$ normalized eigenvalues $\xi_N(\lambda_{ \hat k _{j-1} +i}({\bM}_N) ), \, 1 \leq i \leq k_j$ are eigenvalues of a same matrix of the form $V_{k_j, N} + o_{\mathbb P}(1) $ and then that \begin{equation}{\label{egjoint}}
(\xi_N(\lambda_{\hat k _{j-1} + i}({\bM}_N) ); \, 1 \leq i \leq k_j)= (\lambda_i(V_{k_j, N}); \, 1 \leq i \leq k_j) + o_{\mathbb P}(1).
\end{equation}
\noindent Note that the authors do not develop  this difficulty  in \cite{BY2} (pp. 464-465).
Hence, in the last step of the proof (Step 4 in Section \ref{preuves}), we detail the additional arguments which are needed to get \eqref{egjoint} when $k_j >1$.

 Our approach will cover Cases a) and b) and we will handle both cases once this will be possible. In fact, the main difference appears in the proof of the convergence in distribution of the matrix $V_{k_j, N}$ which gives rise to the "occurrence or non-occurrence" of the distribution $\mu$ in the limiting fluctuations and then justifies the non-universality (resp. universality) in Case a) (resp. b)).\\

The proof is organized in four steps as follows. In Steps 1 and 2, we explain how to obtain (\ref{ecriture1adeX}): we exhibit the matrix ${\check{X}}_{k_{j},N}$ and bring its
 leading term $V_{k_{j},N}$ to light in Step 2. We establish the convergence in distribution of the matrix $V_{k_{j},N}$ in Step 3.
Step 4 is devoted to the concluding arguments of the proof.



\section{Proofs of Theorem \ref{CLT}, Theorem \ref{CLTkN} and Theorem \ref{rank1} }\label{preuves}
As far as possible, we handle  both the  proofs of Theorem \ref{CLT} and Theorem \ref{CLTkN}. We will proceed in four steps. First, let us introduce a few notations.\\
For any matrix $M \in {\cal M}_N(\C)$, we denote by $\Tr$ (resp. $\tr_N$) the classical (resp. normalized) trace. For a rectangular matrix, $||M ||$ is  the operator norm of $M$ and $||M||_{HS}:= (\Tr(MM^*))^{1/2}$ the Hilbert-Schmidt norm. \\
For an Hermitian matrix, we denote by ${\rm Spect}(M)$  the spectrum of $M$. 
For $z \in \C \backslash {\rm{Spect}}(M)$,  $G_{M}(z)~=~(zI_N - M)^{-1} $ denotes the resolvent of $M$ (we suppress the index $M$ when there is no confusion). We have the following :
\begin{equation} \label{lem0}
  \mbox{ For }
x  >\lambda_{1} (M); \qquad \qquad 
\Vert G(x) \Vert \leq \frac{1}{ x- \lambda_{1}(M)}.
\end{equation}
For a $m \times q$ matrix $B$ (or $\bf B$)  and some integers $
1 \leq p \leq m$ and $1 \leq l \leq q$, we denote respectively by $[B]^\nwarrow_{p \times l}$,
$[B]^\nearrow_{p \times l}$, $[B]^\swarrow_{p \times l}$ and $[B]^\searrow_{p \times l}$   the upper left, upper right,
lower left and  lower right
 corner of size $p\times l$  of the matrix $B$. If $p=l$, we will often replace the indices $p \times l$
 by $p$ for convenience. Moreover if $p=m$ , we may replace $\nearrow $ or $\searrow$ by
 $\rightarrow$ and $\swarrow$ or $\nwarrow$ by $\leftarrow$. Similarly if $l = q$,
 we may replace $\nearrow $ or $\nwarrow$ by
 $\uparrow$ and $\swarrow$ or $\searrow$ by $\downarrow$.

 \noindent For simplicity in the writing we will define the $k\times k$, resp. $N-k \times N-k $, resp. $k\times N-k$  matrix $W_k$, resp. $W_{N-k}$, resp. $Y$, by setting
 \begin{equation}{\label{decompoWN}}
\bW_N=
\left( \begin{array}{ll}
W_k ~~Y \\ Y^* ~ W_{N-k}
\end{array}\right).
\end{equation}
\noindent Given ${\bf B} \in \mathcal M _N (\mathbb C)$, we will denote by ${\tilde \bB}$ the $N \times N$ matrix
given by
\begin{equation*} \label{DefTildeB}
\tilde {\bB} := {\rm diag}(U_k^*, I_{N-k}) \, {\bB} \, {\rm diag}(U_k, I_{N-k})
= \left( \begin{array}{ll}
\tilde{B}_k & \tilde{B}_{k\times N-k} \\ \tilde{B}_{N-k \times k} &  \tilde{B}_{N-k}\end{array}\right).
\end{equation*}
One obviously has that $\tilde{B}_{N-k}={B}_{N-k}$.\\
In this way, we define the matrices $\tilde{\bM} _N$, $\tilde{\bW}_N$ and $\tilde{\bA} _N$. In particular, we notice from (\ref{hypA_N}) that
\begin{equation} \label{DefTildeA}
\tilde{\bA}_N=
 {\rm diag}
       (\theta_j I_{k_j}, (\theta_{l}I_{k_l})_{l \leq J_{+ \sigma}, l \not= j}, Z_{N-k_{+ \sigma}}) =  \left( \begin{array}{ll}
\tilde{A}_k & \tilde{A}_{k\times N-k} \\ \tilde{A}_{N-k \times k} &  {A}_{N-k}
\end{array}\right).
\end{equation}
Note also that since ${A}_{N-k}$ is a submatrix of $Z_{N-k_{+ \sigma}}$, all its eigenvalues are strictly smaller than $ \sigma$.\\
Let  $0< \delta < ({\rho_{\theta_j} -2 \sigma})/{2}$. For any random variable $\lambda$, define the events
 \begin{eqnarray*}\Omega_N^{(1)}(\lambda)& =&\left\{ \lambda_1\left(
       \frac{\bW_N}{\sqrt{N}}  + {\rm diag} (U_k, I_{N-k}) \, {\rm diag}
       (0_{k_{+ \sigma}}, Z_{N-k_{+ \sigma}}) \, {\rm diag} (U_k^*, I_{N-k}) \right) <
     2 \sigma+ \delta; \lambda > \rho_{\theta_j} - \delta  \right\},\nonumber \\
 \Omega_N^{(2)}& =&\left\{ \lambda_1\left( \frac{W_{N-k}}{\sqrt{N}}+{A}_{N-k} \right) \leq 2 \sigma+ \delta
  \right\},\end{eqnarray*}
  \noindent and
\begin{equation}{\label{OmegaNLambda}}
\Omega_N(\lambda)=\Omega_N^{(1)}(\lambda)\bigcap
 \Omega_N^{(2)} .
\end{equation}

\noindent On  $\Omega_N(\lambda)$, neither $\lambda$ nor $\rho_{\theta_j}$ are  eigenvalues of $ M_{ N-k
}:=\frac{W_{N-k}}{\sqrt{N}}+{A}_{N-k}$, thus   the resolvent $\widehat{G}(x)$ of $ M_{ N-k}$ is well defined
at $x = \lambda$ and  $x= \rho_{\theta_j}$.

\noindent Note that from Theorem \ref{ThmASCV}, for any random sequence $\Lambda_N$ converging towards $\rho_{\theta_j}$ in probability, $\lim_{N \vers \infty} \mathbb
P(\Omega_N(\Lambda_N)) = 1$.\\

\noindent Let us now introduce on $\Omega_N(\lambda)$ some auxiliary matrices   that will be of basic use to the proofs. 
\begin{equation} \label{defBN} B_{k,N}= W_k +
\frac{1}{\sqrt{N}}\Big(Y\widehat{G}(\rho_{\theta_j})Y^* -(N-k)
\frac{\sigma^2}{\theta_j}I_k\Big).\end{equation}
\begin{equation} \label{defVN}
V_{k_{+\sigma},N}:=\left[U_k^* B_{k,N} U_k\right]^\nwarrow_{k_{+\sigma}}.
\end{equation}
\begin{equation}\label{deftauN}\tau_N(\lambda)=\frac{1}{N}
(\lambda- \rho_{\theta_j})Y \widehat{G}(\lambda) \widehat{G}(\rho_{\theta_j})^2 Y^*, \end{equation}
\begin{equation} \label{defphiN}\phi_N=
- \frac{1}{N}\left( Y\widehat{G}(\rho_{\theta_j})^2 Y^* -   {\sigma^2}
\Tr \widehat{G}(\rho_{\theta_j})^2 I_k \right), \end{equation}
\begin{equation} \label{defpsiN}\psi_N=
 - \sigma^2
\frac{N-k}{N} \Big(\tr_{N-k}\widehat{G}(\rho_{\theta_j})^2 -
\frac{1}{\theta_j^2- \sigma^2} \Big)I_k,\end{equation}
\begin{equation} \label{defDN}
{c_{\theta_j}} D_{k,N}(\lambda)=\tau_N (\lambda)+\phi_N+\psi_N.
\end{equation}
\begin{equation} \label{defTN}T_N (\lambda)=
\left[U_k^* ( W_k  +   \frac{1}{\sqrt{N}}Y \hat{G}(\lambda) Y^*) U_k                             \right]^\nearrow_{k_{+\sigma} \times (k-k_{+\sigma})},\end{equation}
\begin{equation} \label{defDeltaN}\Delta_N(\lambda) = \left[U_k^* Y \hat{G}(\lambda) \tilde{A}_{N-k \times k} \right]^\nearrow_{k_{+\sigma} \times (k-k_{+\sigma})},\end{equation}
 \begin{equation} \label{defGamma}
\Gamma_{k_{+\sigma} \times k-k_{+\sigma}}(\lambda) = T_N(\lambda) + \Delta_N(\lambda).
\end{equation}
\begin{equation}{\label{DefQkN}}
\bQ_{k,N}( \lambda):= \tilde M_k  + \tilde M_{k \times N-k} \hat G (\lambda) \tilde M_{ N-k \times k },
\end{equation}
\begin{eqnarray}\label{DefSigma}
\Sigma_{k-k_{+\sigma}}(\lambda) & = & \left( \left[\bQ_{k,N}(\lambda)\right]^\searrow_{k-k_{+\sigma}}
- \lambda I_{k-k_{+\sigma}}\right)^{-1}.
\end{eqnarray}
Note that we will justify that $ \Sigma_{k-k_{+\sigma}}(\lambda)$ is well defined in the course of the proof of 
Proposition \ref{ecriture0X} below.  Finally, set
\begin{eqnarray}\label{defX}\bX _{k_{+\sigma},N}(\lambda)&=&\left[U_k^* B_{k,N}U_k\right]^\nwarrow_{k_{+\sigma}}+\sqrt{N}
{\rm diag} \left( 0_{k_j}, (\theta_l- \theta_j) I_{k_l}, \, l=1, \ldots,{J_{+ \sigma}}, l \not= j \right) \nonumber \\ &&+ \xi_N(\lambda)
\left[U_k^* D_{k,N}(\lambda)U_k\right]^\nwarrow_{k_{+\sigma}} \nonumber \\
&&+\left( \frac{\sigma^2}{\theta_j^2- \sigma^2}\frac{\xi_N(\lambda)}{c_{\theta_j}}
 \frac{k}{{N}}-\frac{k}{\sqrt{N}} \frac{\sigma^2}{\theta_j}\right)I_{k_{+\sigma}} \nonumber \\
 &&  - \frac{1}{\sqrt{N}}\Gamma_{k_{+\sigma} \times k-k_{+\sigma}}(\lambda)\Sigma_{k-k_{+\sigma}}(\lambda) \Gamma_{k_{+\sigma} \times k-k_{+\sigma}}(\lambda)^*.\end{eqnarray}


\noindent {\bf STEP 1:} We show that an eigenvalue of $\bM_N$ is an eigenvalue of a matrix of size $k_{+\sigma}$. More, precisely, we have:
\begin{proposition} \label{ecriture0X}
For any random variable $\lambda$ and any $k_{+\sigma} \times k_{+\sigma}$ random matrix ${\Delta} _{k_{+\sigma}}$,
on $\Omega_N(\lambda)$, $\lambda$ is an eigenvalue of $\tilde{\bM}_N + {\rm diag}({\Delta} _{k_{+\sigma}}, 0)$
iff $\xi_N(\lambda)$ is an eigenvalue of $ \bX _{k_{+\sigma},N}(\lambda)+\sqrt{N} {\Delta}_{k_{+\sigma}}$ where $\bX _{k_{+\sigma},N}(\lambda)$ is defined by  (\ref{defX}).
Moreover, the $k-k_{+\sigma}\times k-k_{+\sigma}$ matrix $\Sigma_{k-k_{+\sigma}}(\lambda)$ defined by \eqref{DefSigma} is such that $$ \Vert \Sigma_{k-k_{+\sigma}}(\lambda) \Vert \leq {1}/{(\rho_{\theta_j}-2 \sigma -2 \delta)}.$$
\end{proposition}

\noindent
{\bf Proof:} 
Let $\lambda$ be a random variable. On $\Omega_N(\lambda)$,
\begin{eqnarray*}
\det( \bM_N - \lambda I_N) & = & \det( \tilde \bM_N - \lambda I_N)\\
& = &  \det \left( \begin{array}{ll}
\tilde M_k - \lambda I_k & \tilde M_{k \times N-k}\\
M_{ N-k \times k } &  M_{N-k} - \lambda I_{N-k}
\end{array}\right)\\
& = & \det(  M_{N-k} - \lambda I_{N-k}) \det \left( \tilde M_k - \lambda I_k + \tilde M_{k \times N-k} \hat G (\lambda) \tilde M_{ N-k \times k }\right).
\end{eqnarray*}
The last equality in the above equation follows from \eqref{calculdet}.
Since on
$\Omega_N(\lambda)$, $\lambda$ is not an eigenvalue of $ M_{ N-k
}$, we can deduce that $\lambda$ is an eigenvalue of $\tilde{M}_N$
 if and only if it is an eigenvalue of
$$
\bQ_{k,N}( \lambda)= \tilde M_k  + \tilde M_{k \times N-k} \hat G (\lambda) \tilde M_{ N-k \times k }.
$$
Now, note that we have also from \eqref{calculdet} that
\begin{eqnarray*}
& &
 \det \left( \left[  \frac{\tilde {\bf W} _N}{\sqrt N}  \right]^\searrow_{N-k_{+\sigma}} + Z_{N-k_{+\sigma}} - \lambda I_{N-k_{+\sigma}} \right)\\
 & = & \det \left(     \frac{W _{N-k}}{\sqrt N} + \left[ Z_{N-k_{+\sigma}}\right]^\searrow_{N-k} - \lambda I_{N-k}
\right) \times \det \left ( \left[\bQ_{k,N}(\lambda)\right]^\searrow_{k-k_{+\sigma}} - \lambda I_{k-k_{+\sigma}} \right).
\end{eqnarray*}
The matrix $\left[ \frac{\tilde {\bW}_N}{\sqrt N}
\right]^\searrow_{N-k_{+\sigma}} + Z_{N-k_{+\sigma}}$ is a submatrix of $ \frac{\tilde {\bW}_N}{\sqrt
  N} + {\rm diag}(0_{k_{+\sigma}},Z_{N-k_{+\sigma}})$ whose eigenvalues are (on $\Omega_N(\lambda)$) smaller than $2 \sigma + \delta$. So, since on $\Omega_N(\lambda)$, $\lambda$ is greater than
$\rho_{\theta_j} - \delta > 2 \sigma + \delta$, we can conclude that $\lambda$ cannot be an eigenvalue of $\left[  \frac{\tilde {\bW}_N}{\sqrt N}
\right]^\searrow_{N-k_{+\sigma}} + Z_{N-k_{+\sigma}}$, and then  neither of $\left[\bQ_{k,N}(\lambda)\right]^\searrow_{k-k_{+\sigma}}$.
\noindent
Thus, we can define
\begin{eqnarray*}
\Sigma_{k-k_{+\sigma}}(\lambda) & = & \left( \left[\bQ_{k,N}(\lambda)\right]^\searrow_{k-k_{+\sigma}}
- \lambda I_{k-k_{+\sigma}}\right)^{-1}.
\end{eqnarray*}
Moreover on $\Omega_N(\lambda)$, one can see using \eqref{calculdet} that
if $\lambda_0$ is an eigenvalue of $\left[\bQ_{k,N}(\lambda)\right]^\searrow_{k-k_{+\sigma}}
- \lambda I_{k-k_{+\sigma}}$ then $\lambda$ is an eigenvalue of
$$\Big[  \frac{\tilde \bW_N}{\sqrt{N}}   \Big]^\searrow_{N-k_{+\sigma}} + Z_{N-k_{+\sigma}}- {\rm{diag}}(
{\lambda_0} I_{k-k_{+\sigma}},0_{N-k}).$$
Hence,
$$\lambda \leq \lambda_1\Big(\Big[   \frac{\tilde \bW_N}{\sqrt{N}}  \Big]^\searrow_{N-k_{+\sigma}} + Z_{N-k_{+\sigma}} \Big) + {\vert \lambda_0 \vert}$$
and then
$${\vert \lambda_0 \vert} \geq \rho_{\theta_j} - \delta - 2 \sigma - \delta,$$
so that  finally
\begin{equation}\label{borneSigma} \Vert \Sigma_{k-k_{+\sigma}}(\lambda) \Vert \leq \frac{1}{\rho_{\theta_j}-2 \sigma -2 \delta}. \end{equation}
Using oncemore \eqref{calculdet}, we get that on $\Omega_N(\lambda)$, $\lambda$ is an eigenvalue of
$\bQ_{k,N}(\lambda)$ if and only if it is an eigenvalue of
$\left[\bQ_{k,N} (\lambda)\right]^\nwarrow_{k_{+\sigma}}- \left[\bQ_{k,N} (\lambda) \right]^\nearrow_{k_{+\sigma} \times
  k-k_{+\sigma}}\Sigma_{k-k_{+\sigma}}( \lambda)\left[{\bQ_{k,N} (\lambda)}\right]^\swarrow_{ k-k_{+\sigma} \times k_{+\sigma} }$ or equivalently if and only if $\xi_N(\lambda)$  is an eigenvalue of
$$c_{\theta_j}\sqrt{N} \left( \left[\bQ_{k,N}(\lambda)\right]^\nwarrow_{k_{+\sigma}} - \rho_{\theta_j}I_{k_{+\sigma}} - \left[\bQ_{k,N}(\lambda)\right]^\nearrow_{k_{+\sigma} \times
  k-k_{+\sigma}}\Sigma_{k-k_{+\sigma}}( \lambda)\left[\bQ_{k,N}(\lambda)\right]^\swarrow_{k-k_{+\sigma} \times k_{+\sigma} } \right).$$
Now using
    $$ \widehat{G}(\lambda)-\widehat{G}(\rho_{\theta_j})= - (\lambda - \rho_{\theta_j})
    \widehat{G}(\rho_{\theta_j})\widehat{G}(\lambda),$$
 one can replace   $ \widehat{G}(\lambda)$ by $ \widehat{G}(\rho_{\theta_j})+\left[-(\lambda -\rho_{\theta_j})\widehat{G}(\rho_{\theta_j})\left(\widehat{G}(\rho_{\theta_j}) -(\lambda -\rho_{\theta_j})\widehat{G}(\rho_{\theta_j})
\widehat{G}(\lambda) \right)\right]$ and get the following writing
\begin{equation}\label{TS}
\frac{1}{\sqrt{N}}Y\widehat{G}(\lambda)Y^* = \frac{1}{\sqrt{N}}Y\widehat{G}(\rho_{\theta_j})Y^*
+ \xi_N(\lambda)D_{k,N}(\lambda) - \xi_N(\lambda)\frac{N-k}{N} \frac{\sigma^2}{c_{\theta_j}(\theta_j^2- \sigma^2)}I_k \end{equation}
\noindent where
 \begin{eqnarray*}
{c_{\theta_j}} D_{k,N}(\lambda)&=& \frac{1}{N}
(\lambda- \rho_{\theta_j})Y \widehat{G}(\lambda) \widehat{G}(\rho_{\theta_j})^2 Y^*
- \frac{1}{N}\left( Y\widehat{G}(\rho_{\theta_j})^2 Y^* -   {\sigma^2}
\Tr \widehat{G}(\rho_{\theta_j})^2 I_k \right)\\&&
 - \sigma^2
\frac{N-k}{N} \Big(\tr_{N-k}\widehat{G}(\rho_{\theta_j})^2 -
\frac{1}{\theta_j^2- \sigma^2} \Big)I_k.
\end{eqnarray*}
\noindent Then
\begin{eqnarray*}c_{\theta_j}\sqrt{N} \left ( \left[\bQ_{k,N} (\lambda) \right]^\nwarrow_{k_{+\sigma}} - \rho_{\theta_j}I_{k_{+\sigma}} \right) &=&c_{\theta_j}\left\{
\left[U_k^* \left(W_{k}+ \frac{1}{\sqrt{N}}\Big(Y\widehat{G}(\rho_{\theta_j})Y^* -(N-k)
\frac{\sigma^2}{\theta_j}I_k\Big) \right)U_k\right]^\nwarrow_{k_{+\sigma}}\right.\\&&
 +\sqrt{N}
{\rm diag} \left( 0_{k_j}, (\theta_l- \theta_j) I_{k_l}, \, l=1, \ldots,{J_{+ \sigma}}, l \not= j \right)\\
&&\left. + \xi_N(\lambda)
\left[U_k^* D_{k,N}(\lambda)U_k\right]^\nwarrow_{k_{+\sigma}}  -\frac{k}{\sqrt{N}} \frac{\sigma^2}{\theta_j}I_{k_{+\sigma}}  +\frac{\sigma^2}{\theta_j^2- \sigma^2}\frac{\xi_N(\lambda)}{c_{\theta_j}}
 \frac{k}{{N}}I_{k_{+\sigma}}\right\}\\
&&-\frac{\sigma^2}{\theta_j^2- \sigma^2}{\xi_N(\lambda)}I_{k_{+\sigma}}.
 \end{eqnarray*}
The  proposition (adding an extra matrix $\Delta_{k_{+\sigma}}$ for future computations) readily follows. $\Box$


\vspace{.3cm}
Throughout Steps 2 and 3, $\Lambda_N$ denotes any random sequence converging in probability towards $\rho_{\theta_j}$.
The aim of these two steps is to study the limiting behavior of the matrix $\bX _{k_{+\sigma},N}(\Lambda_N)$ (defined by (\ref{defX})) as $N$ goes to infinity.\\

\noindent
\noindent {\bf STEP 2:}
We first focus on the negligible terms in $\bX _{k_{+\sigma},N}(\Lambda_N)$ and establish the following.
\begin{proposition}\label{ecriture2deX} Assume that $k \ll \sqrt N$. For any random sequence $\Lambda_N $ converging in probability towards $\rho_{\theta_j}$, on $\Omega_N(\Lambda_N)$,
\begin{equation}\label{estimX}
\bX _{k_{+\sigma},N}(\Lambda_N) = V_{k_{+\sigma},N} + \sqrt{N}
{\rm diag} \left( 0_{k_j}, (\theta_l- \theta_j) I_{k_l}, \, l=1, \ldots,{J_{+ \sigma}}, l\not= j \right) + (1+|\xi_N(\Lambda _N)|)^2 o_{\mathbb P}(1),
\end{equation}
with $V_{k_{+\sigma},N}$ 
defined by (\ref{defVN})
\end{proposition}

The proof of this proposition is quite long and is divided in several lemmas.
Although our final result in the case $k$ infinite holds only for $k \ll \sqrt N$, we will give some estimates for $k \ll N$ once this is possible.

\begin{lemme}\label{etudedeDNkN} Let $k \ll N$. Then,  on $\Omega_N(\Lambda_N)$,
\begin{equation} \label{NeglD}
\left[U_k^* D_{k,N}(\Lambda_N) U_k\right]^\nwarrow_{k_{+\sigma}}
=o_{\mathbb P}(1).
\end{equation}
 \end{lemme}

\noindent {\bf Proof of Lemma \ref{etudedeDNkN}:} 
$D_{k,N}(\Lambda_N), \ \tau_N, \ \phi_N$ and $\psi_N$ are respectively defined by (\ref{defDN}), (\ref{deftauN}), (\ref{defphiN}) and (\ref{defpsiN}) .

\noindent
Since $Y$ is a submatrix of $W_N$, $\limsup_N \frac{\Vert Y \Vert }{\sqrt{N}} \leq 2$ and  $\frac{\Vert Y \Vert }{\sqrt{N}} = O_{\mathbb P}(1)$. Thus,
\begin{equation} \label{normetauN}
 \Vert \tau_N \Vert = O_{\mathbb P}(\Lambda_N - \rho_{\theta_j}) = o_{\mathbb P}(1).
 \end{equation}
 where we used that $\Vert \hat{G}(\lambda)\Vert \leq  \frac{1}{(\rho_{\theta_j} - 2 \sigma- 2 \delta)}$ for $\lambda = \rho_{\theta_j}$ or $\lambda = \Lambda_N$.
Therefore,
$$ \Vert \left[U_k^* \tau_N U_k\right]^\nwarrow_{k_{+\sigma}} \Vert =o_{\mathbb P}(1).$$

\noindent It follows from Lemma \ref{LemmaAppendix} in the Appendix that $$\left[U_k^* \psi_N U_k \right]^\nwarrow_{k_{+\sigma}}:=  - \sigma^2
\frac{N-k}{N} \Big[\tr_{N-k}\widehat{G}(\rho_{\theta_j})^2 -
\frac{1}{\theta_j^2- \sigma^2} \Big]I_{k_{+\sigma}} =o_{\mathbb P}(1).$$

\noindent Now, we have
\begin{eqnarray*}
\mathbb{E}\left(\Vert \left[U_k^* \phi_N U_k\right]^\nwarrow_{k_{+\sigma}} \1 \Vert^2_{HS} \right)
& \leq & \mathbb{E}\left(\Vert \left[U_k^* \phi_N U_k\right]^\nwarrow_{k_{+\sigma}} \deux  \Vert^2_{HS} \right)\\
 &\leq& \sum_{p,q=1}^{k_{+ \sigma}} \frac{1}{N^2}\mathbb{E}\left(\vert {\cal U}(p)^*\widehat{G}(\rho_{\theta_j})^2{\cal U}(q) - \sigma^2 \Tr
 \widehat{G}(\rho_{\theta_j})^2 \delta_{p,q} \vert^2 \deux   \right),
 \end{eqnarray*}
 \noindent where for any $p=1,\ldots,k_{+ \sigma}$, we let $ {\cal U}(p)=~^t[ (Y^*U_k)_{1, p}, \ldots, (Y^*U_k)_{N-k, p}]$.
 We first state some properties of the vectors ${\cal U}(p)$.

\begin{lemme} \label{defcalU}
Let ${\cal U}$ denote the $  N-k \times k_{+\sigma}$ matrix $[Y^*U_k ]^\leftarrow_{k_{+\sigma}}$. Then, the rows $({\cal U}_{i.} ; i \leq N-k)$ are centered i.i.d vectors in $\C^{k_{+\sigma}}$, with a distribution depending on $N$. Moreover, we have for all $1 \leq p,q \leq k_{+\sigma}$:
\begin{eqnarray}
&& \E({\cal U}_{1p} {\bar{\cal U}}_{1q}) =\delta_{p,q} \sigma^2 \quad \text{with} \quad \E({\cal U}_{1p} {\cal U}_{1q}) =  0 \, \text{ in the complex case,} \nonumber \\&&
\E[|{\cal U}_{ip}|^2 |{\cal U}_{iq}|^2] = (1+ \frac{t}{2}\delta_{p,q} )\sigma^4 +  [ \E(|W_{12}|^4) - (1+ \frac{t}{2}) \sigma^4] \sum_{l=1}^k |(U_k)_{l,p}|^2 |(U_k)_{l,q}|^2. \label{pU}
\end{eqnarray}
Since $\sum_{l=1}^k |(U_k)_{l,p}|^4 \leq 1$, the fourth moment of ${\cal U}_{1p}$ is uniformly bounded.
\end{lemme}
\noindent We skip the proof of this lemma which follows from straightforward computations using the independence of the entries of $Y$ and the fact that $U_k$ is unitary. \\
Then, according to Theorem \ref{BaiSilver98} and using  \eqref{lem0},
\begin{eqnarray*} \frac{1}{N^2}
\mathbb{E}\left (\vert  {\cal U}(p)^*\widehat{G}(\rho_{\theta_j})^2{\cal U}(p) -
\sigma^2 \Tr
 \widehat{G}(\rho_{\theta_j})^2  \vert^2 \1  \right)&\leq &\frac{K}{N} \E\left(  \tr_N  \widehat{G}(\rho_{\theta_j})^4 \deux  \right) \\& \leq & \frac{K}{N} \E\left(\Vert \widehat{G}(\rho_{\theta_j})\Vert^4
\deux  \right)\\
& \leq & \frac{K}{N} \, \frac{1}{(\rho_{\theta_j}-2 \sigma -
\delta)^4}.
\end{eqnarray*}
Besides for $p \neq q$, using the independence between $({\cal U}(p),{\cal U}(q))$ and $\widehat{G}(\rho_{\theta_j})$, we have:
\begin{eqnarray*}
\mathbb{E}\left(\vert {\cal U}(p)^*\widehat{G}(\rho_{\theta_j})^2{\cal U}(q)  \vert^2 \deux   \right)
&= & \sum_{i,j,l,m}^{N-k} \E[ \bar{{\cal U}}_{ip} (G^2)_{ij} {\cal U}_{jq}{\cal U}_{lp}\overline{ (G^2)_{lm}}\bar{{\cal U}}_{mq}\deux ] \\
 &= & \sum_{i,j,l,m}^{N-k} \E[ \bar{{\cal U}}_{ip}  {\cal U}_{jq}{\cal U}_{lp}{{\bar{\cal U}}}_{mq}] \E[ (G^2)_{ij}\overline{(G^2)_{lm}}\deux ]
\end{eqnarray*}
where we denote by $G$ the matrix $\widehat{G}(\rho_{\theta_j})$ for simplicity.
>From Lemma \ref{defcalU}, for $p \neq q$, the only terms giving a non null expectation in the above equation are those for which:
\begin{itemize}
\item[1)] $i = l $, $j=m$ and $i\neq j$. In this case,
$$  \E[ \bar{{\cal U}}_{ip}  {\cal U}_{jq}{\cal U}_{ip}\bar{{\cal U}}_{jq}]   = \E[ \bar{{\cal U}}_{ip} {\cal U}_{ip}]  \E[ {\cal U}_{jq}\bar{{\cal U}}_{jq}] = \sigma^4$$
and
$$\sum_{i,j,i\neq j}^{N-k} \E[ (G^2)_{ij}\overline{(G ^2)_{ij}}\deux ] \leq \E \Tr( G^4\deux ).$$
\item[2)] $i=j=k=l$. In this case, using (\ref{pU}), there is a constant $C>0$ such that
$$  \E[ \bar{{\cal U}}_{ip}  {\cal U}_{iq}{\cal U}_{ip}\bar{{\cal U}}_{iq}]   = \E[|{\cal U}_{ip}|^2 |{\cal U}_{iq}|^2]
 \leq C.$$
Moreover $$\sum_{i=1}^{N-k} \E[ G^2_{ii}{\bar G^2}_{ii}\deux ] \leq \E \Tr (G^4 \deux ),$$
 \end{itemize}
 Therefore,
\begin{equation} \label{Maj} \mathbb{E}\left(\vert {\cal U}(p)^*\widehat{G}(\rho_{\theta_j})^2{\cal U}(q)  \vert^2 \deux   \right)
\leq (C + \sigma^4)  \ \E  \Tr (\widehat{G}(\rho_{\theta_j})^4\deux ).\end{equation}
 Hence,$$\frac{1}{N^2} \mathbb{E}\left(\vert {\cal U}(p)^*\widehat{G}(\rho_{\theta_j})^2{\cal U}(q)  \vert^2 \deux   \right)
\leq  \frac{ C+ \sigma^4}{N} \E\left(\Vert \widehat{G}(\rho_{\theta_j})\Vert^4
\deux  \right) \leq  \frac{C+ \sigma^4}{N} \, \frac{1}{(\rho_{\theta_j}- 2 \sigma -
\delta)^4}.$$
Thus \begin{eqnarray*}
\mathbb{E}\left( \Vert \left[U_k^* \phi_N U_k\right]^\nwarrow_{k_{+\sigma}} \Vert^2 \1 \right)
 & \leq &(C+\sigma^4)\frac{k_{+ \sigma}^2}{N} \, \frac{1}{(\rho_{\theta_j}- 2 \sigma -
\delta)^4} \label{HSphi}.\end{eqnarray*}
The convergence in probability of  $\left[U_k^* \phi_N U_k\right]^\nwarrow_{k_{+\sigma}}$ towards zero  readily follows by Tchebychev inequality. \\
Lemma \ref{etudedeDNkN} is established. $\Box$

\paragraph{} For simplicity, we now write
$$\Sigma(\Lambda_N)=\Sigma_{k-k_{+\sigma}}(\Lambda_N).$$
Let us define
\begin{equation}\label{Defreste}
R_{k,N}(\Lambda_N): =-  \frac{k}{\sqrt{N}} \frac{\sigma^2}{\theta_j}I_{k_{+\sigma}} +\frac{\sigma^2}{ \theta_j^2 - \sigma^2} \frac{\xi_N(\Lambda_N)}{c_{\theta_j}}  \frac{k}{{N}} I_{k_{+\sigma}}  - \frac{1}{\sqrt{N}} \Gamma_{k_{+\sigma} \times
  k-k_{+\sigma}}(\Lambda_N) \Sigma(\Lambda_N) \Gamma_{k_{+\sigma} \times k-k_{+\sigma}}(\Lambda_N)^*.
\end{equation}
To get Proposition \ref{ecriture2deX}, it remains to prove that if $k \ll \sqrt{N}$,
\begin{equation}\label{reste}
R_{k,N}(\Lambda_N)=(1+ \vert \xi_N(\Lambda_N)\vert )^2 o_{\mathbb P}(1).
\end{equation}

Once $k \ll \sqrt{N}$, we readily have that
$$ -  \frac{k}{\sqrt{N}} \frac{\sigma^2}{\theta_j}I_{k_{+\sigma}}+\frac{\sigma^2}{ \theta_j^2 - \sigma^2} \frac{\xi_N(\Lambda_N)}{c_{\theta_j}}  \frac{k}{{N}} I_{k_{+\sigma}}= (1+ \vert \xi_N(\Lambda_N)\vert )^2 o_{\mathbb P}(1).$$
\noindent Hence, (\ref{reste})  will follow if we
prove
\begin{lemme} \label{borneQ2}  Assume that $k \ll \sqrt{N}$. Let $ \Gamma_{k_{+\sigma} \times
  k-k_{+\sigma}}(\lambda)$ and $\Sigma(\lambda)$ be defined as  (\ref{defGamma}) and (\ref{DefSigma}).
  On $\Omega_N(\Lambda_N)$,
\begin{equation}\label{Q2a} \frac{1}{\sqrt{N}} \Gamma_{k_{+\sigma} \times
  k-k_{+\sigma}}(\Lambda_N) \Sigma(\Lambda_N) \Gamma_{k_{+\sigma} \times k-k_{+\sigma}}(\Lambda_N)^*
= (1+ \vert \xi_N(\Lambda_N) \vert )^2 o_{\mathbb P}(1).
\end{equation}

\end{lemme}

\noindent
For the proof, we use the following decomposition ($T_N(\lambda)$ and $\Delta_N(\lambda)$ being defined by (\ref{defTN}) and
(\ref{defDeltaN})):
\begin{eqnarray} \label{decompoQ}
\lefteqn{\Gamma_{k_{+\sigma} \times k-k_{+\sigma}}(\Lambda_N)\Sigma (\Lambda_N)  \Gamma_{{k_{+\sigma} \times k-k_{+\sigma}}}(\Lambda_N)^*} \nonumber \\
&& =T_N \Sigma T_N^*
+T_N \Sigma \Delta_N(\Lambda_N)^* + \Delta_N(\Lambda_N) \Sigma \Delta_N(\Lambda_N)^* + \Delta_N(\Lambda_N) \Sigma T_N^*
\end{eqnarray}
where (using (\ref{TS}))
\begin{eqnarray*}
T_N:=T_N (\Lambda_N) &=&
\left[U_k^* ( W_k  +   \frac{1}{\sqrt{N}}Y \hat{G}(\Lambda_N) Y^* ) U_k \right]^\nearrow_{k_{+\sigma} \times (k-k_{+\sigma})} \\
&=& [U_k^* B_{k,N}U_k]^\nearrow_{k_{+\sigma} \times k-k_{+\sigma}} + \xi_N(\Lambda_N)
\left[U_k^*  D_{k,N}(\Lambda_N) U_k\right]^\nearrow_{k_{+\sigma} \times k-k_{+\sigma}}
\end{eqnarray*}
and we replaced $\Sigma (\Lambda_N) $ by $\Sigma$.
We will prove the following lemma on $T_N$.

\begin{lemme}\label{etudedeTN} 
If $k \ll N$, 
\begin{equation} \label{lemma3.8}
\Vert\left[U_k^*  D_{k,N}(\Lambda_N) U_k\right]^\nearrow_{k_{+\sigma} \times k-k_{+\sigma}}\Vert=  o_{\mathbb P}(1).
\end{equation}
If $ k \ll \sqrt{N}$,
\begin{equation}
\label{estimBrectangle}
\Vert\left[U_k^*  B_{k,N}(\Lambda_N) U_k\right]^\nearrow_{k_{+\sigma} \times k-k_{+\sigma}}\Vert_{HS}= o_{\mathbb P}(N^{\frac{1}{4}}).
\end{equation}
and therefore, for $k \ll \sqrt{N}$,
$$\Vert T_N \Vert = o_{\mathbb P}(N^{\frac{1}{4}})(1+ \vert \xi_N(\Lambda_N)\vert ).$$
\end{lemme}

\noindent {\bf Proof of Lemma \ref{etudedeTN}:}
To prove \eqref{lemma3.8}, we use the decomposition $${c_{\theta_j}} \left[U_k^*  D_{k,N}(\Lambda_N) U_k \right]^\nearrow_{k_{+\sigma} \times k-k_{+\sigma}}=\left[U_k^* \tau_NU_k\right]^\nearrow_{k_{+\sigma} \times k-k_{+\sigma}}+ \left[U_k^*\phi_NU_k\right]^\nearrow_{k_{+\sigma} \times k-k_{+\sigma}}.$$
As in the proof of Lemma \ref{etudedeDNkN}, we have
\begin{eqnarray*}
\mathbb{E}\left( \Vert \left[U_k^* \phi_N U_k\right]^\nearrow_{k_{+\sigma} \times k-k_{+\sigma}}\Vert_{HS}^2 \1 \right)&\leq &
(C+\sigma^4) \frac{kk_{+\sigma}}{N} \, \frac{1}{(\rho_{\theta_j}- 2 \sigma -
\delta)^4} ,\end{eqnarray*}
\noindent so that, for  $k \ll N$ and using Tchebychev inequality,
we can deduce that $$ \Vert \left[U_k^* \phi_N U_k\right]^\nearrow_{k_{+\sigma} \times k-k_{+\sigma}}\Vert_{HS}\1= o_{\mathbb P}(1).$$
>From \eqref{normetauN},
$$  \Vert\left[U_k^*  \tau_N U_k\right]^\nearrow_{k_{+\sigma} \times k-k_{+\sigma}}\Vert = o_{\mathbb P}(1).$$
and therefore
 $$ \Vert\left[U_k^*  D_{k,N}(\Lambda_N) U_k\right]^\nearrow_{k_{+\sigma} \times k-k_{+\sigma}}\Vert= 
 o_{\mathbb P}(1).$$
Thus, \eqref{lemma3.8} is established.\\
For \eqref{estimBrectangle}, recall that $
[U_k^* B_{k,N}U_k]^\nearrow_{k_{+\sigma} \times k-k_{+\sigma}}= [U_k^* W_k U_k]^\nearrow_{k_{+\sigma} \times k-k_{+\sigma}}+
\frac{1}{\sqrt{N}}[U_k^*
Y\widehat{G}(\rho_{\theta_j})Y^*U_k]^\nearrow_{k_{+\sigma} \times
  k-k_{+\sigma}}$. 
  Since
$\Vert W_k\Vert = O_{\mathbb P}(\sqrt{k})$, we have $\Vert [U_k^* W_k U_k]^\nearrow_{k_{+\sigma} \times k-k_{+\sigma}}\Vert^2 =  O_{\mathbb P}(\sqrt{k})$. 
Hence, as $k \ll \sqrt{N}$, we can deduce that $\Vert [U_k^* W_k U_k]^\nearrow_{k_{+\sigma} \times k-k_{+\sigma}}\Vert_{HS}= o_{\mathbb P}(N^{\frac{1}{4}})$.\\
Now, let us prove the same estimate for the remaining term. Using  the same proof as in \eqref{Maj}, one can get that for $p \neq q$, for some constant $C>0$, $$\mathbb{E}\left(\vert {\cal U}(p)^*\widehat{G}(\rho_{\theta_j}){\cal U}(q)  \vert^2 \deux   \right)\leq C \E  \Tr (\widehat{G}(\rho_{\theta_j})^2\deux )$$
and then  that for some constant $C>0$,
$$ \E[\Vert \frac{1}{\sqrt{N}}[U_k^* Y\widehat{G}(\rho_{\theta_j})Y^*U_k]^\nearrow_{k_{+\sigma} \times k-k_{+\sigma}}\Vert_{HS}^2 \deux] \leq C k k_{+ \sigma} \frac{1}{ (\rho_{\theta_j}-2\sigma -\delta)^2}. $$
Then using that
\begin{equation*} \label{normeHS}
\PP \left(\Vert \frac{1}{\sqrt{N}}[U_k^* Y\widehat{G}(\rho_{\theta_j})Y^*U_k]^\nearrow_{k_{+\sigma} \times k-k_{+\sigma}}\Vert_{HS} \deux > \epsilon N^{\frac{1}{4}}  \right)
\leq \frac{1}{\epsilon^2 \sqrt{N}} \E[\Vert \frac{1}{\sqrt{N}}[U_k^* Y\widehat{G}(\rho_{\theta_j})Y^*U_k]^\nearrow_{k_{+\sigma} \times k-k_{+\sigma}}\Vert_{HS}^2]
\end{equation*}
we deduce since $k \ll \sqrt{N}$ that
$$\Vert \frac{1}{\sqrt{N}}[U_k^* Y\widehat{G}(\rho_{\theta_j})Y^*U_k]^\nearrow_{k_{+\sigma} \times k-k_{+\sigma}}\Vert_{HS} \deux= o_{\mathbb P}(N^{\frac{1}{4}}).$$
Thus \eqref{estimBrectangle} and Lemma \ref{etudedeTN} are proved. $\Box$  \\

Using that
\begin{equation} \label{SigmaN}
\Vert \Sigma \Vert \leq \frac{1}{\rho_{\theta_j}-2 \sigma -2 \delta},
\end{equation}
one can readily notice that Lemma \ref{etudedeTN} leads to  \begin{equation} \label{TN} \frac{1}{\sqrt{N}}T_N \Sigma T_N^* = (1+
  \vert \xi_N(\Lambda_N)\vert )^2o_{\mathbb P}(1). \end{equation}
We now consider the remaining terms in the r.h.s of
  $(\ref{decompoQ})$. We first show the following result where we recall that $\Delta_N(\rho_{\theta_j}) = [U_k^* Y
\hat{G}(\rho_{\theta_j}) \tilde{A}_{N-k \times
  k}]^\nearrow_{k_{+\sigma} \times k-k_{+\sigma}}$.

\begin{lemme}\label{avecrho}
$\frac{1}{\sqrt{N}}T_N \Sigma \Delta_N(\rho_{\theta_j})^*$, $\frac{1}{\sqrt{N}}\Delta_N(\rho_{\theta_j}) \Sigma \Delta_N(\rho_{\theta_j})^*$ and $\frac{1}{\sqrt{N}}\Delta_N(\rho_{\theta_j}) \Sigma T_N^*$ are all equal to some $ (1+ \vert \xi_N(\Lambda_N)\vert )o_{\mathbb P}(1).$
\end{lemme}

\noindent {\bf Proof of Lemma \ref{avecrho} :} We will show that, on $\Omega_N(\Lambda_N)$, for any $u >0$,
\begin{equation}\label{deltak1N}\Delta_N(\rho_{\theta_j})=o_{\mathbb P}(N^{u}).\end{equation}
One can readily see that this leads to the announced result combining Lemma \ref{etudedeTN}, (\ref{SigmaN}) and (\ref{deltak1N}). \\
First, using the fact that $U_k^* Y$ is independent of $\deux \hat{G}(\rho_{\theta_j})$ and that
for any $p$, the random vector $ {\cal U}(p)=^t[ (Y^*U_k)_{1,p},
\ldots, (Y^*U_k)_{N-k, p}]$ has independent centered entries with variance $\sigma^2$, one has that
\begin{eqnarray*}\mathbb{E}(\1 \Tr
  \Delta_N(\rho_{\theta_j})\Delta_N(\rho_{\theta_j})^*) &\leq&
\mathbb{E}(\deux \Tr \Delta_N(\rho_{\theta_j})\Delta_N(\rho_{\theta_j})^*)
\\ &= &k_{{+\sigma}} \sigma^2  \mathbb{E}\left\{\deux \Tr   [\hat{G}^2(\rho_{\theta_j}) \tilde{A}_{N-k \times k-k_{+\sigma}} \tilde{A}_{N-k \times k- k_{+\sigma}}^*]
\right\}\\
&\leq& k_{{+\sigma}} \sigma^2
\mathbb{E} \left\{ \deux  \Vert \hat{G}(\rho_{\theta_j}) \Vert^2  \Tr
\tilde{A}_{N-k \times k-k_{+\sigma}} \tilde{A}_{N-k \times k-k_{+\sigma}}^* \right \}\\
& \leq & \frac{k_{{+\sigma}} \sigma^2}{(\rho_{\theta_j}-2 \sigma - \delta)^2}  \Tr
{A}_{N}^2  \\
& = &  \frac{k_{{+\sigma}} \sigma^2}{(\rho_{\theta_j}-2 \sigma - \delta)^2} \sum_{l=1}^J k_l \theta_l^2.
\end{eqnarray*}
Therefore, $\mathbb{P}( \1 \Vert  \Delta_N(\rho_{\theta_j}) \Vert_{HS} \geq \epsilon N^{u})
\leq {\epsilon^{-2} N^{-2u}} \mathbb{E}(\1 \Vert  \Delta_N(\rho_{\theta_j})\Vert_{HS}^2)$ goes to zero as $N$ tends to infinity. Hence (\ref{deltak1N}) holds true on $\Omega_N(\Lambda_N)$ and the proof of Lemma
\ref{avecrho} is complete. $\Box$\\

Let us now prove that
\begin{lemme}\label{diff}
$\Delta_N(\Lambda_N)=\Delta_N(\rho_{\theta_j})+
 O_{\mathbb P}(\vert \xi_N(\Lambda_N)\vert ) $.
\end{lemme}

\noindent {\bf Proof of Lemma \ref{diff}:} We have
$$\Delta_N(\Lambda_N)-\Delta_N(\rho_{\theta_j})=
-(\Lambda_N -\rho_{\theta_j}) [U_k^* Y
\hat{G}(\rho_{\theta_j}) \hat{G}(\Lambda_N) \tilde{A}_{N-k \times
  k}]^\nearrow_{k_{+\sigma} \times k-k_{+\sigma}} . $$
Let us define $\nabla_{k_{+\sigma}}= [U_k^* Y
\hat{G}(\rho_{\theta_j}) \hat{G}(\Lambda_N) \tilde{A}_{N-k \times
  k}]^\nearrow_{k_{+\sigma} \times k-k_{+\sigma}}$. Then for some constant $C>0$ depending on the matrix $\tilde{A}_{N-k \times
  k}$,
 \begin{eqnarray*}
 \Tr( \nabla_{k_{+\sigma}}\nabla_{k_{+\sigma}}^*) &\leq &C \Vert \hat{G}(\rho_{\theta_j})\Vert^2
\Vert \hat{G}(\Lambda_N)\Vert ^2 \Tr({\cal U}^*{\cal U}) \\
&\leq & \frac{C}{(\rho_{\theta_j}-2 \sigma - 2\delta)^4} \Tr({\cal U}^*{\cal U})
\end{eqnarray*}
where we denote as before ${\cal U}=[Y^*U_k ]^\leftarrow_{k_{+\sigma}}$. Thus letting $C':= C \, c_{\theta_j}^{-2}$,
$$ \Vert \Delta_N(\Lambda_N)-\Delta_N(\rho_{\theta_j})\Vert_{HS}^2 \leq C' (\xi_N(\Lambda_N))^2  \frac{1}{(\rho_{\theta_j}-2 \sigma - 2\delta)^4} \frac{1}{N} \Tr({\cal U}^*{\cal U}).$$
>From Lemma \ref{defcalU}, it  follows  that
 $$\frac{1}{N}  \Tr({\cal U}^*{\cal U})  \stackrel{\mathbb P}{\rightarrow} k_{+\sigma} \sigma^2$$
implying that $ \Vert  \Delta_N(\Lambda_N)-\Delta_N(\rho_{\theta_j})\Vert_{HS} =
 O_{\PP} (\vert \xi_N(\Lambda_N)\vert ). \; \Box$

\paragraph{} We are now in position to conclude the proof of Lemma \ref{borneQ2}. Indeed, writing
$$\Delta_N(\Lambda_N)\Sigma T_N^* = \left(\Delta_N(\Delta_N)-\Delta_N(\rho_{\theta_j}) \right) \Sigma T_N^* + \Delta_N(\rho_{\theta_j}) \Sigma T_N^*$$
and \begin{eqnarray*} \Delta_N(\Lambda_N)\Sigma \Delta_N(\Lambda_N)&=&\Delta_N(\rho_{\theta_j})\Sigma \Delta_N(\rho_{\theta_j})^*\\
&& +\left(\Delta_N(\Lambda_N)-\Delta_N(\rho_{\theta_j})\right)\Sigma \Delta_N(\rho_{\theta_j})^*\\
&&+ \left(\Delta_N(\Lambda_N)-\Delta_N(\rho_{\theta_j}) \right) \Sigma \left(\Delta_N(\Lambda_N)-\Delta_N(\rho_{\theta_j}) \right)^*\\
&&+\Delta_N(\rho_{\theta_j})\Sigma \left(\Delta_N(\Lambda_N)-\Delta_N(\rho_{\theta_j}) \right)^*,
\end{eqnarray*}
we deduce from Lemmas \ref{etudedeTN}, \ref{diff} and \eqref{SigmaN},
\eqref{deltak1N} that
$\frac{1}{\sqrt{N}} \Delta_N(\Lambda_N)\Sigma T_N^*$ and $\frac{1}{\sqrt{N}}\Delta_N(\Lambda_N)\Sigma \Delta_N(\Lambda_N)^*$
are both equal to some $(1+ \vert \xi_N(\Lambda_N)\vert )o_{\mathbb P}(1)$.
Using also (\ref{TN}), we can deduce that
\begin{equation}  \frac{1}{\sqrt{N}}  \Gamma_{k_{+\sigma}\times k-k_{+\sigma}}(\Lambda_N)\Sigma
\Gamma_{k_{+\sigma}\times k-k_{+\sigma}}(\Lambda_N)^*= (1 + \vert \xi_N(\Lambda_N)\vert )^2
  o_{\mathbb P}(1)
\end{equation}
which gives (\ref{Q2a}) and completes the proof of Lemma \ref{borneQ2}. $\Box$\\

Combining all the preceding, we have established Proposition \ref{ecriture2deX}. We now prove that provided it converges in distribution, with a probability going to one as $N$ goes to infinity, $\xi_N(\Lambda_N)$ is actually an eigenvalue of a matrix of size $k_j$.
\begin{lemme}\label{normeV} For all $u>0$,
$$\frac{\Vert V_{k_{+\sigma}, N} \Vert_{HS}}{{N}^u}=o_{\mathbb P}(1).$$
\end{lemme}
\noindent {\bf Proof:} Straightforward computations lead to the existence of some constant $C$ such that
$$\E \left( \Vert \left[U_k^*W_{k}U_k\right]_{k_{+\sigma}} \Vert_{HS} \right) \leq C.$$
The convergence of ${\Vert \left[U_k^*W_{k}U_k\right]_{k_{+\sigma}} \Vert}/{{N}^u}$ in probability towards zero readily follows by Tchebychev inequality.
Following the proof in Lemma \ref{etudedeDNkN} of the convergence in probability of $\left[U_k^*\Phi_{N}U_k\right]_{k_{+\sigma}}$ towards zero, one can get that
\begin{eqnarray*}
\mathbb{E}\left( \Vert \left[U_k^* \frac{1}{\sqrt{N}}\deux \Big(Y\widehat{G}(\rho_{\theta_j})Y^*-\sigma^2 \Tr\widehat{G}(\rho_{\theta_j}) I_k\Big) U_k\right]^\nwarrow_{k_{+\sigma}} \Vert^2 \1 \right)
 & \leq &\frac{(C+\sigma^4)k_{+ \sigma}^2} {(\rho_{\theta_j}- 2 \sigma -
\delta)^2},\end{eqnarray*}
\noindent and the convergence in probability towards zero of the term inside the above expectation follows by Tchebychev inequality.
Since moreover according to Lemma \ref{LemmaAppendix},$$\frac{1}{\sqrt{N}}\deux \left( \Tr\widehat{G}(\rho_{\theta_j})-(N-k) \frac{1}{\theta_j} \right) = o_{\mathbb P}(1),$$
\noindent we can deduce that $$N^{-u}{\Vert \left[U_k^* \frac{1}{\sqrt{N}}\deux \Big(Y\widehat{G}(\rho_{\theta_j})Y^*- (N-k) \frac{\sigma^2}{\theta} I_k\Big) U_k\right]^\nwarrow_{k_{+\sigma}} \Vert}\deux =o_{\mathbb P}(1).$$
\noindent The proof of Lemma \ref{normeV} is complete. $\Box$

\begin{proposition}\label{ecriture3deX} Let $\Delta_{k_j}$ be an arbitrary $k_j \times k_j$ random matrix.  If $\xi_N(\Lambda_N)$ converges in distribution, then, with a probability going to one as $N$ goes to infinity, it  is an eigenvalue of $\bX_{k_{+\sigma},N}(\Lambda_N)  + {\rm diag}( \Delta_{k_j},0) $ iff $\xi_N(\Lambda_N)$ is an eigenvalue of a matrix  $\check{X}_{k_j, N}(\Lambda_N)  + \Delta_{k_j}$ of size $k_j$, satisfying
\begin{equation}{\label{estXter}}
\check{X}_{k_j, N} (\Lambda_N)  = V_{k_j, N} +  o_{\mathbb P}(1)
\end{equation}
where $V_{k_j,N}$ is the $k_j \times k_j$ element in the block decomposition of  $V_{k_{+\sigma}, N}$ defined by \eqref{defVN}; namely
$$V_{k_j,N} = U_{K_j \times k_j}^* [B_{k,N}]^\nwarrow _{K_j}U_{K_j \times k_j} $$
\noindent with $U_{K_j \times k_j}$ and $B_{k,N}$ defined respectively by \eqref{defUj} and \eqref{defBN}.
\end{proposition}

\noindent {\bf Proof of Proposition \ref{ecriture3deX}:}

Since $\xi_N(\Lambda_N)$ converges in distribution, we can write the matrix $\bX_{k_{+\sigma},N}(\Lambda_N)$ given by \eqref{estimX} as
$$ \bX_{k_{+\sigma},N} (\Lambda_N)  = \sqrt{N} {\rm diag}(
0_{k_j}, ((\theta_l- \theta_j)I_{k_{l}})_{l \not= j}) + \check{R}_{k_{+\sigma},N}(\Lambda_N)$$
where $\check{R}_{k_{+\sigma},N}(\Lambda_N) := V_{k_{+\sigma},N} + o_{\mathbb P}(1)$.
Let us decompose $\bX_{k_{+\sigma},N}(\Lambda_N) $ in blocks as
$$\bX_{k_{+\sigma},N} (\Lambda_N) = \left( \begin{array}{ll}
 X_{k_j,N}   & X_{k_{j} \times k_{+\sigma}-k_j,N } \\
 X_{k_{+\sigma}-k_j \times k_j,N } & X_{k_{+\sigma}-k_j,N}
\end{array}\right).$$
We first show that $\xi_N(\Lambda_N)$ is not an eigenvalue of $X_{k_{+\sigma}-k_j,N  } $. Let $\alpha = \inf_{l\not= j} | \theta_l - \theta_j| >0$. Since,
$$ X_{k_{+\sigma}-k_j,N  } =  \sqrt{N} {\rm diag}(
((\theta_l- \theta_j)I_{k_{l}})_{l\not= j}) + \check{R}_{k_{+\sigma}-k_j,N },$$
if $\mu$ is an eigenvalue of $X_{k_{+\sigma}-k_j }$, then
$${|\mu|} /{\sqrt{N}} \geq  \alpha - || \check{R}_{k_{+\sigma}-k_j ,N}|| / \sqrt{N}.$$
Now, using Lemma \ref{normeV},
$$||\check{R}_{k_{+\sigma}-k_j,N }|| / \sqrt{N} = o_{\mathbb P}(1).$$
Hence $\xi_N(\Lambda_N)$ cannot be an eigenvalue of
$X_{k_{+\sigma}-k_j,N } $.
Therefore, we can define
\begin{eqnarray*}
\check{X}_{k_j, N} &= & X_{k_j,N}  - X_{k_{j} \times k_{+\sigma}-k_j,N }
(X_{k_{+\sigma}-k_j ,N } - \xi_N(\Lambda_N) I_{k_{+\sigma}-k_j  })^{-1} X_{k_{+\sigma}-k_j \times k_j,N }\\
 &= & V_{k_j, N}   - \check{R}_{k_{j} \times k_{+\sigma}-k_j,N }
(X_{k_{+\sigma}-k_j ,N } - \xi_N(\Lambda_N)  I_{k_{+\sigma}-k_j  })^{-1}
\check{R}_{k_{+\sigma}-k_j \times k_j,N }+o_{\mathbb P}(1).
\end{eqnarray*}
To get \eqref{estXter}, it remains to show that
$$ ||  \check{R}_{k_{j} \times k_{+\sigma}-k_j,N }
(X_{k_{+\sigma}-k_j ,N } - \xi_N(\Lambda_N)  I_{k_{+\sigma}-k_j  })^{-1}
\check{R}_{k_{+\sigma}-k_j \times k_j,N } || = o_{\mathbb P}(1).$$
This follows from the previous computations showing that (for some constant $C>0$)
$$ || (X_{k_{+\sigma}-k_j,N  } - \xi_N(\Lambda_N)  I_{k_{+\sigma}-k_j  })^{-1} || \leq \left({C +o_{\mathbb P}(1)}\right)/{\sqrt{N}},$$
combined with the definition of $\check{R}_{k_{+\sigma},N}(\Lambda_N)$ and Lemma \ref{normeV}.
The statement of the proposition then follows from \eqref{calculdet}. $\Box$ \\

\noindent {\bf STEP 3:}
We now examine the convergence of the $k_{j} \times k_{j}$ matrix $V_{k_j,N}=U_{K_j \times k_j}^*[B_{k,N}]^\nwarrow _{K_j}U_{K_j \times k_j}$

\begin{proposition}\label{convBinfinite}
 The $k_{j}\times k_{j}$ matrix
$V_{k_j,N}=U_{K_j \times k_j}^*[B_{k,N}]^\nwarrow _{K_j}U_{K_j \times k_j}$ converges in distribution to a GU(O)E($k_{j}\times k_{j},  \frac{\theta_j^2 \sigma^2}{\theta_j^2- \sigma^2})$ if and only if $\max_{p=1}^{k_j}\max_{i=1}^{K_j} \vert (U_k)_{ip}\vert $
converges to zero when $N$ goes to infinity.
\end{proposition}
\noindent{\bf Proof}\\
Assume that $\max_{p=1}^{k_j}\max_{i=1}^{K_j} \vert (U_k)_{ip}\vert $
converges to zero when $N$ goes to infinity. We decompose the proof of the convergence 
of $U_{K_j \times k_j}^*[B_{k,N}]^\nwarrow _{K_j}U_{K_j \times k_j}$ in distribution to a GU(O)E($k_{j}\times k_{j},  \frac{\theta_j^2 \sigma^2}{\theta_j^2- \sigma^2})$ into the two following lemmas.

\begin{lemme}\label{premierbout}
If $\max_{p=1}^{k_j}\max_{i=1}^{K_j} \vert (U_k)_{ip}\vert $ converges to zero when $N$ goes to infinity  then the $k_{j}\times k_{j}$ matrix
$U_{K_j \times k_j}^*[W_{k}]^\nwarrow _{K_j}U_{K_j \times k_j}$ converges in distribution to a GU(O)E($k_{j}\times k_{j},  \sigma^2 )$.
\end{lemme}
\noindent{\bf Proof of Lemma \ref{premierbout}:} First we consider the complex case. Let $\alpha_{pq} \in \C$, $1 \leq p<q \leq
k_{j}$ and $\alpha_{pp} \in \R$, $1 \leq p \leq
k_{j}$, and define
$$ L_N (\alpha) := \sum_{1\leq p<q  \leq {k}_{j}} (\alpha_{pq} (U_k^* W_k U_k)_{pq} + \overline{\alpha_{pq} (U_k^* W_k U_k)_{pq}}) + \sum_{1\leq p \leq {k}_{j}} 2\alpha_{pp} (U_k^* W_k U_k)_{pp}.$$
\noindent We have
$$  L_N (\alpha) =\sum_{i=1}^{K_j} D_i (W_N)_{ii} + \sum_{1\leq i<l \leq {K_j}} R_{il} (\sqrt{2}\Re e
((W_N)_{il})) +\sum_{1\leq i<l \leq {K_j}} I_{il} (\sqrt{2}\Im m
((W_N)_{il})),$$
\noindent where
$$D_i
=2\Re e \bigg ( \sum_{1\leq p\leq q \leq {k}_{j}} \alpha_{pq} (U_k)_{iq} \overline{(U_k)}_{ip} \bigg),$$
$$R_{il}= \sqrt{2} \Re e \bigg( \sum_{1\leq p\leq q \leq {k}_{j}} \alpha_{pq} ((U_k)_{lq} \overline{(U_k)}_{ip} + (U_k)_{iq} \overline{(U_k)}_{lp})\bigg),$$
$$I_{il}= \sqrt{2} \Im m \bigg( \sum_{1\leq p\leq q \leq {k}_{j}} \overline{\alpha_{pq}} (\overline{(U_k)}_{lq} {(U_k)_{ip}} - \overline{(U_k)}_{iq} {(U_k)_{lp}})\bigg).$$
\noindent Hence $ L_N (\alpha)=\sum_{m=1}^{{K_j}^2}\beta_{m,N} \phi_m$ where $\phi_m$ are i.i.d random variables with distribution $\mu$ and $\beta_{m,N}$ are real constants (depending on the $\alpha_{pq}$) which satisfy $\max_{m=1}^{{K_j}^2} |\beta_{m,N}| \rightarrow 0$
when $N$ goes to infinity.  Therefore the cumulants of $L_N (\alpha)$ are given by
 $C_n^{(N)}= \sum_{m=1}^{{K_j}^2} \beta_{m,N}^n C_n(\mu)$ for any $n \in \mathbb{N}^*$ where $C_n(\mu)$ denotes the $n$-th cumulant of $\mu$ (all are  finite since $\mu$ has moments of any order).
 In particular $C_1^{(N)}=0$.  We are going to  prove that the variance of
$L_N (\alpha)$ is actually constant, given by
\begin{equation}\label{varianceL}\frac{C_2^{(N)}}{\sigma^2} = \sum_{m=1}^{K_j^2} \beta_{m,N}^2   = 2
\sum_{1\leq p< q \leq {k}_{j}} \vert \alpha_{pq} \vert^2 +4 \sum_{1\leq p\leq {k}_{j}} \vert \alpha_{pp} \vert^2.\end{equation}
 One may  rewrite $L_N (\alpha)$ as $$L_N (\alpha)= \Tr \left(H U_k^* W_k U_k \right)$$
 \noindent where $H$ is the $k\times k$ Hermitian matrix defined by
 $$H_{pq}=\alpha_{qp} ~~\mbox{~~if~~} p>q  ~~\mbox{~~and~~} H_{pp}=2\alpha_{pp}.$$
 \noindent Hence using 
 $$\E \left[ (W_k)_{ji} (W_k)_{qp}\right]= \delta_{jp} \delta_{iq} \sigma^2$$
 \noindent it is easy to see that 
\begin{eqnarray*}\E \left[ L_N(\alpha)^2 \right]&=&\sigma^2  \Tr \left[\left(U_k H U_k^*\right)^2\right]\\&=&
\sigma^2 \Tr H^2.
\end{eqnarray*}
Then (\ref{varianceL}) readily follows.
In the following, we let $const=\sum_{m=1}^{{K_j}^2}\beta_{m,N}^2$.\\
 Since $\vert C_n^{(N)}\vert \leq  const \max_{m=1}^{K_j^2} \vert \beta_{m,N}\vert ^{n-2} \vert C_n(\mu)\vert $,
$C_n^{(N)}$ converges to zero for each $n \geq 3$.
Thus we can
deduce from Janson's theorem \cite{J} that
  $L_N (\alpha)$ converges
 to a centered gaussian distribution with variance $\sigma^2 ( 2
\sum_{1\leq p< q \leq {k}_{j}} \vert \alpha_{pq} \vert^2 +4 \sum_{1\leq p\leq {k}_{j}} \vert \alpha_{pp} \vert^2) $ and the proof of Lemma \ref{premierbout} is complete in the complex case.\\

 \noindent Dealing with symmetric matrices, one needs to consider the random variable
$$ L_N (\alpha) := \sum_{1\leq p <q \leq {k}_{j}} \alpha_{pq} (U_k^t W_k U_k)_{pq}  + \sum_{1\leq p \leq {k}_{j}} \alpha_{pp} (U_k^t W_k U_k)_{pp}$$
for any real numbers $\alpha_{pq}, \, p \leq q$. One can similarly prove that $ L_N (\alpha)$ converges
 to a centered gaussian distribution with variance $\sigma^2 ( 2
\sum_{1\leq p< q \leq {k}_{j}}  \alpha_{pq}^2 +2 \sum_{1\leq p\leq {k}_{j}}  \alpha_{pp} ^2). $ $\Box$
\begin{remarque} Note that Lemma \ref{premierbout} is true under the assumption of the existence of a fourth moment. This can be shown by using a Taylor development of the Fourier transform of $L_N( \alpha)$.
\end{remarque}

\begin{lemme}\label{deuxiemebout}
If $\max_{p=1}^{{k}_{j}}\max_{i=1}^{K_j} \vert (U_k)_{ip}\vert $ converges to zero when $N$ goes to infinity then the $k_{j}\times~k_{j}$ matrix
$\frac{1}{\sqrt{N}}U_{K_j\times k_j}^*\left[\Big(Y\widehat{G}(\rho_{\theta_j})Y^* -(N-k)
\frac{\sigma^2}{\theta_j}I_k\Big)\right]^{\nwarrow}_{K_j}U_{K_j\times k_j}$ converges towards a
GU(O)E($k_{j}\times k_{j}, \frac{\sigma^4}{\theta_j^2-\sigma^2})$.
\end{lemme}
\vspace{.3cm}
\noindent
{\bf Proof of Lemma \ref{deuxiemebout}:} We shall apply a slightly modified version of Theorem \ref{Silver} (see Theorem 7.1 in \cite{BY2}) but requiring the finiteness of sixth (instead of fourth) moments. Let $K = k_{j}(k_{j}+1)/{2}$. The set $\{1, \ldots K\}$ is indexed by $l= (p,q)$ with $1 \leq p \leq q \leq {k}_{j}$, taking the lexicographic order. We define a sequence of i.i.d centered vectors $(x_i, y_i)_{i \leq N-k}$ in $\C^K \times \C^K$ by $x_{li} = {\cal U}_{ip}$ and $y_{li} = {\cal U}_{iq}$ for $l= (p, q)$ where ${\cal U}$ is defined in Lemma \ref{defcalU}. The matrix $A$ of size $N-k$ is the matrix $\widehat{G}(\rho_{\theta_j})$ and is independent of ${\cal U}$.
Note that we are not exactly in the context of Theorem 7.1 of \cite{BY2} since the i.i.d vectors $(x_i, y_i)_i$ depend on $N$ (and should be rather denoted by $(x_{i, N}, y_{i,N})_i$) but their distribution satisfies:
\begin{enumerate}
\item $\rho(l) = \E[\bar{x}_{l1}y_{l1}] = \delta_{p,q} \sigma^2$ for $l=(p,q)$ is independent of $N$.
\item $\E[\bar{x}_{l1}y_{l'1}] = \delta_{p,q'} \sigma^2$ if $l=(p,q), l'= (p', q')$ (see $B_2$ in \eqref{defBmatrix}).
\item Complex case: $\E[\bar{x}_{l1}\bar{x}_{l'1}] =\E[y_{l1}y_{l'1}] = 0$ if $l=(p,q), l'= (p', q')$ (see $B_3$ in \eqref{defBmatrix}).

Real case: $\E[\bar{x}_{l1}\bar{x}_{l'1}] = \sigma^2 \delta_{p,p'}$ and $\E[y_{l1}y_{l'1}] = \sigma^2 \delta_{q,q'}$ if $l=(p,q), l'= (p', q')$.
\item  (see $B_1$ in \eqref{defBmatrix})
$$\left\{ \begin{array}{lcl}
\E[\bar{x}_{l1}y_{l1}\bar{x}_{l'1}y_{l'1}] &= &\sigma^4(\delta_{p,q}\delta_{p',q'} + \delta_{p,q'}\delta_{p',q}) + \\
& & [ \E(|W_{12}|^4) - 2\sigma^4] \sum_{i=1}^{K_j} (U_k)_{i,q}\overline{(U_k)_{i,p}}(U_k)_{i,q^{'}}\overline{(U_k)_{i,p^{'}}} \;
 \mbox{ in the complex case, } \\ &&\\
\E[\bar{x}_{l1}y_{l1}\bar{x}_{l'1}y_{l'1}] &= &\sigma^4(\delta_{p,q}\delta_{p',q'} + \delta_{p,q'}\delta_{p',q}+
\delta_{p,p'}\delta_{q,q'} ) + \\
& & [ \E(|W_{12}|^4) - 3\sigma^4] \sum_{i=1}^{K_j} (U_k)_{i,q}\overline{(U_k)_{i,p}}(U_k)_{i,q^{'}}\overline{(U_k)_{i,p^{'}}} \;
 \mbox{ in the real case. }
 \end{array} \right. $$
 \end{enumerate}
 Under the assumption that $\max_{p=1}^{{k}_{j}}\max_{i=1}^{K_j} \vert (U_k)_{i,p}\vert $
converges to zero when $k$ goes to infinity, the last term in the r.h.s of the two above equations tends to 0.

\noindent
It can be seen that the proof of Theorem 7.1 still holds in this case once we verify that for $\epsilon >0$ and for $z = x$ or $y$, for any $l$,
\begin{equation} \label{moment4uniforme}
\E[|z_{l1}|^4 1\! \!{\sf I}_{(|z_{l1}| \geq \epsilon N^{1/4})}] \longrightarrow 0 \quad \text{as } {N \rightarrow \infty}.
\end{equation}
We postpone the proof of \eqref{moment4uniforme} to the end of the proof. Assuming that \eqref{moment4uniforme} holds true, we obtain the CLT theorem 7.1 (\cite{BY2}): the Hermitian matrix $Z_N=(Z_N(p,q))$ of size $k_{j}$ defined by
$$Z_N(p,q) =  \frac{1}{\sqrt{N-k}} [\sum_{i, i' = 1, \ldots N-k} \bar{{\cal U}}_{ip} \widehat{G}(\rho_{\theta_j})_{ii'} {\cal U}_{i'q} - \delta_{p,q} \sigma^2  \Tr(\widehat{G}(\rho_{\theta_j}))]$$
converges to  an Hermitian  Gaussian matrix $G$. The Laplace transform of $G$ (considered as a vector of $\C^K$, that is of $\{G_{pq}, \, 1 \leq p \leq q \leq {k}_{j}\}$) is given for any $c \in \mathbb C ^K$ by
$$\E[ \exp(c^T G)] = \exp[ \frac{1}{2} c^T B c]$$ where the $K \times K$ matrix
$B=(B(l,l'))$ is given by:
$B= \lim_N B_1(N) + B_2 + B_3$ with
\begin{eqnarray*}
B_1(N)(l,l') &=& \omega ( \E[\bar{x}_{l1}y_{l1}\bar{x}_{l'1}y_{l'1}] - \rho(l) \rho(l')), \\
B_2(l,l') &=& (\theta -\omega)\E[\bar{x}_{l1}y_{l'1}] \E[\bar{x}_{l'1}y_{l1}]  \\
B_3(l,l') &=& (\tau - \omega) \E[\bar{x}_{l1}\bar{x}_{l'1}] \E[y_{l1}y_{l'1}]
\end{eqnarray*}
and the coefficients $\omega, \theta, \tau$ are defined in Theorem \ref{Silver}. Here $A = \widehat{G}(\rho_{\theta_j})$ so that $\omega = {1}/{\theta^2_j}$ and $\theta ={1}/{(\theta^2_j - \sigma^2)}$ (see the Appendix). \\
>From Lemma \ref{defcalU},
$$ B_2(l,l') = (\theta - \omega) \sigma^4 \delta_{p,q'}\delta_{p',q} =  (\theta - \omega) \sigma^4 1_{p=q=p'=q'}.$$
Moreover in the complex case, $B_3 \equiv 0$ and in the real case,
$$B_3(l,l') = (\theta - \omega) \sigma^4 \delta_{l,l'}.$$
>From 4., in the real case,
$$\lim_{N \rightarrow \infty} B_1(N) (l,l') =  \delta_{l,l'} \omega \sigma^4(1+ \delta_{p,q}),$$
and in the complex case,
$$\lim_{N \rightarrow \infty} B_1(N) (l,l') =  \delta_{l,l'} \omega \sigma^4\delta_{p,q}.$$
It follows that $B$ is a diagonal matrix given by:
$$ \left\{ \begin{array} {ll}
B(l,l) = (1+ \delta_{p,q}) \theta \sigma^4 = (1+ \delta_{p,q})\frac{ \sigma^4}{ \theta^2_j - \sigma^2} & \mbox{ in the real case,} \\
B(l,l) =  \delta_{p,q} \theta \sigma^4 = \delta_{p,q}\frac{ \sigma^4}{ \theta^2_j - \sigma^2} & \mbox{ in the complex case.}
\end{array} \right. $$
In the real case, the matrix $B$ is exactly the covariance of the limiting Gaussian distribution $G$. It follows that $G$ is the distribution of the GOE($k_{+\sigma}\times k_{+\sigma}, {\sigma^4}/{(\theta^2_j -\sigma^2)})$.\\
In the complex case, from the form of $B$, we can conclude that the coordinates of $G$ are independent ($B$ diagonal), $G_{pp}$ has variance ${ \sigma^4}/{( \theta^2_j - \sigma^2)}$ and for $p \not= q$, $\Re e(G_{pq})$ and $\Im m(G_{pq})$ are independent with the same variance (since $B(l,l) = 0$ for $p\not= q$).
It remains to compute the variance of $\Re e(G_{pq})$. Since the Laplace transform of $\Re e(Z_N(p,q) )$ and $\Im m(Z_N(p,q) )$ can be expressed as a Laplace transform of $Z_N(p,q) $ and $\overline{Z_N(p,q)} $, we shall apply Theorem 7.1 to $(Z_N(p,q) , \overline{{Z}_N(p,q)} )$ that is to the vectors $x_i = ({\cal U}_{ip}, {\cal U}_{iq})$ and $y_i = ({\cal U}_{iq}, {\cal U}_{ip})$ in $\C^2$. We denote by $\tilde{B}$ the associated "covariance" matrix of size 2. The variance of $\Re e(G_{pq})$ is given by $\frac{1}{2} \lim_{N \rightarrow \infty} \tilde{B}_{12}$ (since
$\tilde{B}_{11} = \tilde{B}_{22} = 0$ from the previous computations) with
$$ \tilde{B}_{12} =  \tilde{B}_{12}(1) +  \tilde{B}_{12}(2) +  \tilde{B}_{12}(3)$$ where here $\tilde{B}_{12}(3) = 0$,
$$ \tilde{B}_{12}(1) =\omega \E[ |{{\cal U}}_{1p}|^2 |{{\cal U}}_{1q}|^2 ] \rightarrow  \omega \sigma^4 \quad \text{and} \quad
  \tilde{B}_{12}(2) =(\theta-\omega) \E[|{\cal U}_{1p}|^2]  \E[|{\cal U}_{1q}|^2] = (\theta-\omega) \sigma^4. $$
Therefore, $var(\Re e(G_{pq})) =\theta \sigma^4 /2= {\sigma^4}/{(2(\theta^2_j -\sigma^2))}$.

 We thus obtain Lemma \ref{deuxiemebout} by using that  $ \Tr(\widehat{G}(\rho_{\theta_j})) = (N-k)  \tr_{N-k}(\widehat{G}(\rho_{\theta_j}))$ and $ \tr_{N-k}(\widehat{G}(\rho_{\theta_j})) \rightarrow {1}/{\theta_j}$. \\
 It remains to prove \eqref{moment4uniforme}. The variable $\alpha_N:= |z_{l1}|^4 1\! \!{\sf I}_{(|z_{l1}| \geq \epsilon N^{1/4})}$ tends to 0 in probability. It is thus enough to prove uniform integrability of the sequence $\alpha_N$, a sufficient condition is given by $\sup_N \E[ \alpha_N^{6/4}] < \infty$. It is easy to see that for any
 $1 \leq p \leq {k}_{j}$, $\sup_N \E[ |{\cal U}_{1p}|^6]  < \infty$ since the Wigner matrix $\bW _N$ has finite sixth moment and $U_k$ is unitary. This proves \eqref{moment4uniforme} and finishes the proof of Lemma \ref{deuxiemebout}. $\Box$\\
~

 Assume now that the matrix $V_{k_j,N}=U_{K_j \times k_j}^*[B_{k,N}]^\nwarrow _{K_j}U_{K_j \times k_j}$ converges in distribution towards a GU(O)E($k_{j} \times k_{j},  \frac{\theta_j^2 \sigma^2}{\theta_j^2- \sigma^2})$ whereas $\max_{p=1}^{k_j}\max_{i=1}^{K_j} \vert (U_k)_{ip}\vert $
does not converge  to zero when $N$ goes to infinity. There exists $p_0 \in \{1,\ldots,k_j\}$ such that 
$\max_{i=1}^{K_j} \vert (U_k)_{ip}\vert$ does not converge  to zero.
Let $i_N$ be such that $\max_{i=1}^{K_j} \vert (U_k)_{ip_0}\vert=\vert (U_k)_{i_N p_0}\vert$.
Now we have
$$\left( V_{k_j,N} \right)_{p_0 p_0}= \vert (U_k)_{i_N p_0}\vert^2 W_{i_N i_N} 
+ X_N$$
\noindent where $X_N$ is a random variable which is independent with 
$\vert (U_k)_{i_N p_0}\vert^2 W_{i_N i_N} $.
One can find a subsequence such that $\vert (U_k)_{i_{\phi(N)} p_0}\vert^2 W_{i_{\phi(N)} i_{\phi(N)}}$
converges in distribution towards $c\xi$ where $c>0$ and $\xi$ is $\mu$-distributed.
This leads to a contradiction using Cramer-L\'evy's Theorem since  $\left( V_{k_j,\phi(N)} \right)_{p_0 p_0}$ converges towards a gaussian variable.
 The proof of Proposition \ref{convBinfinite} is complete.
 $\Box$

\vspace{.3cm}
In the case {\bf a)}, condition of Proposition \ref{convBinfinite} are obviously not satisfied and we have the following asymptotic result.
\begin{proposition}\label{convBfinite}
 In case a), the Hermitian (resp. symmetric) matrix $V _{k_j,N}$ converges in distribution towards the law of $V_{k_j\times k_j}$ of size $k_j$
defined in the following way. Let
$W_{\tilde{K}_j} $ be a Wigner matrix of size ${\tilde{K}_j} $ with distribution given by $\mu$ (cf {\bf (i)}) and $H_{\tilde{K}_j} $
 be a centered Hermitian (resp. symmetric) Gaussian matrix of size ${\tilde{K}_j} $ independent of $W_{\tilde{K}_j} $ with independent entries $H_{pl}$, $p \leq l$ with variance
\begin{equation}\label{variance}
\left\{  \begin{array}{l} \displaystyle v_{pp} = E(H_{pp}^2) = \frac{t}{4}\Big ( \frac{m_4 - 3 \sigma^4}{\theta_j^2} \Big )
   + \frac{t}{2}\frac{\sigma^4 }{\theta_j^2-\sigma^2} ,  \, p= 1, \ldots, {\tilde{K}_j} , \\
\displaystyle v_{pl} = \mathbb{E}(|H_{pl}|^2) =
\frac{\sigma^4}{\theta_j^2-\sigma^2},  \, 1 \leq p<l \leq {\tilde{K}_j} .
\end{array} \right.  \end{equation}
Then, $V_{k_j\times k_j}$
is the  $k_j \times k_j$ matrix defined by
\begin{equation}\label{defV} V_{k_j\times k_j} =  \tilde{U}_{\tilde{K}_j\times k_j}^* (W_{\tilde{K}_j}  + H_{\tilde{K}_j} ) \tilde{U}_{\tilde{K}_j\times k_j} .\end{equation}
\end{proposition}
\noindent The proof follows from Theorem \ref{Silver} and is omitted  since
we have detailed the similar  proof of  Lemma \ref{deuxiemebout}.

\vspace*{0.5cm}
\noindent {\bf STEP 4:} We are now in position to prove that
\begin{equation}{\label{egconv}}
\left ( \xi_N( \lambda_{\hat k _{j-1} + 1} ({\bM}_N) ),
  \ldots, \xi_N( \lambda_{\hat k _{j-1}+ {k_j}} ({\bM}_N))
\right ) \stackrel{\mathcal L }{\longrightarrow} \left(\lambda_1(V_{k_j\times k_j}), \ldots,\lambda_{k_j}(V_{k_j\times k_j})\right).
\end{equation}

To prove (\ref{egconv}), our strategy will be indirect: we start from the matrix $V_{k_j, N}$ and its eigenvalues $(\lambda_i(V_{k_j, N}); \, 1 \leq i \leq k_j)$ and we will reverse the previous reasoning to raise to the normalized eigenvalues $\xi_N(\lambda_{\hat k _{j-1} + i}({\bM}_N) ), \, 1 \leq i \leq k_j$. This approach works in both Cases a) and b) as we now explain.\\

First, for any $1 \leq i \leq k_j$, we define $\Lambda_N^{(i)}$ such that
$$\xi_{N}(\Lambda_N^{(i)})=\lambda_i(V_{k_j,N}),$$
that is $ \Lambda_N^{(i)}= \rho_{\theta_j}+ \lambda_i(V_{k_j,N}) /
c_{\theta_j} \sqrt{N}$. \\
Since $V_{k_j,N}$ converges in distribution towards
$V_{k_j\times k_j}$, $\lambda_i(V_{k_j,N})$ also converges in distribution towards
$\lambda_i(V_{k_j\times k_j})$. Hence $ \xi_{N}(\Lambda_N^{(i)})$ converges in distribution and $\Lambda_N^{(i)}$ converges in probability towards $\rho_{\theta_j}$. Let $
\check{X}_{k_j}^{(i)}\equiv \check{X}_{k_j, N}(\Lambda_N^{(i)}) = V_{k_j, N} +  o_{\mathbb P}(1) $ as defined in Proposition \ref{ecriture3deX}.
This fit choice of $\Lambda_N^{(i)}$ gives that
$$\lambda_i(\check{X}_{k_j}^{(i)})= \xi_N(\Lambda_N^{(i)}) + \epsilon_i, \mbox{~~~~
with $\epsilon_i=o_{\mathbb P}(1)$}.$$
\noindent Hence, $\xi_N(\Lambda_N^{(i)})$ is an eigenvalue of $\tilde{X}_{k_j}^{(i)}-\epsilon_i I_{k_j}$.\\
According to Propositions \ref{ecriture0X} and \ref{ecriture3deX}, on an event $\check{\Omega}_N$ whose probability goes to one as $N$ goes to infinity, there exists some $l_i $ such that
$$\Lambda_N^{(i)}= \lambda_{l_i}\Big({\bM}_N - \frac{\epsilon_i}{\sqrt{N}}\mbox{diag}(I_{k_j},0_{N-{k}_{j}})\Big).$$ The following lines hold on $\check{\Omega}_N$.
By using Weyl's inequalities (Lemma \ref{Weyl}), one has for all $i \in \{ 1, \ldots, k_j \}$ that
$$ \left \vert \lambda_{l_i}\Big({\bM}_N - \frac{\epsilon_i}{\sqrt{N}}\mbox{diag}(I_{k_j},0_{N-{k}_j}) \Big )-\lambda_{l_i}({\bM}_N ) \right \vert \leq \frac{|\epsilon_i|}{\sqrt{N}}.$$
We then deduce that
\begin{equation}{\label{egalEigen}}
 \left ( \xi_N( \lambda_{l_1} ({\bM}_N) ), \ldots, \xi_N( \lambda_{l_{k_j}} ({\bM}_N)) \right )=
\Big(\lambda_1(V_{k_j,N}), \ldots,\lambda_{k_j}(V_{k_j,N})\Big) + o_{\mathbb P}(1)
\end{equation}
and thus
\begin{equation}{\label{cvegalEigen}}
 \left ( \xi_N( \lambda_{l_1} ({\bM}_N) ), \ldots, \xi_N( \lambda_{l_{k_j}} ({\bM}_N)) \right ) \stackrel{\mathcal L }{\longrightarrow}
\Big(\lambda_1(V_{k_j\times k_j}), \ldots,\lambda_{k_j}(V_{k_j \times k_j})\Big).
\end{equation}

Now, to get ({\ref{egconv}}), it is sufficient to prove that
\begin{equation}{\label{probali}}
\mathbb P\left( l_i= \hat k _{j-1}+i; \, i=1, \ldots, k_j \right) \rightarrow 1, \mbox{~~as $N \to \infty$}.
\end{equation}
Indeed, one can notice that on the event $\{  l_i= \hat k _{j-1}+i; \, i=1, \ldots, k_j \}$ the following equality holds true
\begin{equation}{\label{egalbisEigen}}
\left ( \xi_N( \lambda_{\hat k _{j-1} + 1} ({\bM}_N) ),
  \ldots, \xi_N( \lambda_{\hat k _{j-1}+ {k_j}} ({\bM}_N))
\right )= \left ( \xi_N( \lambda_{l_1} ({\bM}_N) ), \ldots,
  \xi_N( \lambda_{l_{k_j}} ({\bM}_N)) \right ).
\end{equation}
Hence, if
$(\ref{probali})$ is satisfied then (\ref{egalbisEigen}) combined with (\ref{cvegalEigen}) imply (\ref{egconv}).\\
\noindent We turn now to the proof of $(\ref{probali})$. The key point is to notice that the $k_j$ eigenvalues of $V_{k_j \times k_j}$
have a joint density. This fact is well-known if $V_{k_j \times k_j}$ is a matrix from the GU(O)E and so when $K_j$ is infinite (Case b)). When $K_j$ is bounded (Case a)), we call on the following
arguments. One can decompose the matrix $\tilde{U}_{\tilde{K}_j \times k_j}^* H_{\tilde{K}_j} \tilde{U}_{\tilde{K}_j \times k_j}$ appearing in
the definition $(\ref{defV})$ of $V_{k_{j}\times k_j}$ in the following way
$$\tilde{U}_{\tilde{K}_j \times k_j}^* H_{\tilde{K}_j} \tilde{U}_{\tilde{K}_j \times k_j}=Q_{k_{j}}+\check{H}_{k_{j}}$$
with $\check{H}_{k_{j}}$ distributed as GU(O)E (using the fact that $\tilde{U}_{\tilde{K}_j \times k_j}^*\tilde{U}_{\tilde{K}_j \times k_j}=I_{k_j}$) and $Q_{k_{j}}$  independent from $\check{H}_{k_{j}}$.
Hence, the law of $V_{k_j\times k_{j}}$
is that of the sum of two random independent matrices: the first one being the matrix $\check{H}_{k_{j}}$ distributed as GU(O)E associated to a Gaussian measure with some variance $\tau$ and the second one being a matrix $Z_{k_j}$ of the form $\tilde{U}_{\tilde{K}_j \times k_j}^* W_{\tilde{K}_j} \tilde{U}_{\tilde{K}_j \times k_j} + Q_{k_{j}}$.
Using the density of the GU(O)E matrix $\check{H}_{k_{j}}$ with respect to  the Lebesgue measure $dM$ on Hermitian (resp. symmetric) matrices, decomposing $dM$ on $\mathbb{U}_N\times (\mathbb{R}^N)_{\leq}$ (denoting by $\mathbb{U}_N $ the unitary (resp. orthogonal) group),
one can easily see that
the distribution of the eigenvalues of $\check{H}_{k_{j}}+Z_{k_j}$ is absolutely continuous with respect to the
Lebesgue measure $d\lambda$ on $\mathbb{R}^n$ with a density given by:
$$f(\lambda_1,\ldots, \lambda_N)= \exp({-\frac{N}{\tau t}\sum_{i=1}^N \lambda_i^2}) \prod_{i<j}(\lambda_i -\lambda_j)^{ \frac{4}{t}}
\mathbb{E}\left(\exp\left\{-\frac{N}{\tau t}\Tr Z_{k_j}^2\right\} I((\lambda_1,\ldots, \lambda_N), Z_{k_j})\right)d\lambda$$
where $I((\lambda_1,\ldots, \lambda_N), Z_{k_j})= \int \exp\left(\frac{2}{\tau t} N\Tr (U \mbox{diag}(\lambda_1,\ldots, \lambda_N)U^*
Z_{k_j})\right) m(dU)$ denoting by $m$ the Haar measure on the unitary (resp. orthogonal) group.\\
Thus, we
deduce that the $k_j$ eigenvalues of $V_{k_j\times k_{j}}$ are distinct (with
probability one). Using Portmanteau's Lemma with (\ref{cvegalEigen}) then implies that the event
$$ \check{\Omega}^{'}_N:= \Big \{  \xi_N( \lambda_{l_1} ({\bM}_N)
)> \cdots > \xi_N( \lambda_{l_{k_j}} ({\bM}_N)) \Big \} \bigcap \check{\Omega}_N$$
is such that $\lim _N \mathbb P ( \check{\Omega}^{'}_N)=1$. By Theorem \ref{ThmASCV}, we notice that the event
$$ \tilde{\Omega}'_N:= \left \{  \lambda_{\hat{k}_{j-1} }
  ({\bM}_N) > \rho _{\theta _j} + \delta > \lambda_{l_1}
  ({\bM}_N) \right \} \bigcap \check{\Omega}^{'}_N \bigcap \left \{
  \lambda_{l_{k_j}} ({\bM}_N) > \rho _{\theta _j} - \delta >
  \lambda_{\hat{k}_{j-1} + k_j+1} ({\bM}_N)  \right \}$$ also
satisfies $\lim _N \mathbb P ( \tilde{\Omega}' _N)=1$, for $\delta$
small enough. This leads to $(\ref{probali})$ since $ \tilde{\Omega}'
_N \subset  \{  l_i= i + \hat{k}_{j-1}, \,  i=1, \ldots, k_j \}.$ \\
The proof of Theorems \ref{CLT} and \ref{CLTkN} is complete. $\Box$ \\
~

 According to Theorem  \ref{CLTkN}, in order to establish Theorem \ref{rank1}, we only need to prove that the condition (\ref{cns}) is actually necessary for universality of the fluctuations. Hence assume that \\$\sqrt{N} (\lambda_{k_1+\cdots+k_{j-1}+1}(\bM_N) - \rho_{\theta_j}) \overset{\mathcal
  L}{\longrightarrow}  {\cal
  N}(0, \frac{t}{2}\sigma^2_{\theta_j})$.  Proposition \ref{ecriture0X}  and Proposition \ref{ecriture3deX} lead to
  $$c_{\theta_j}\sqrt{N} (\lambda_{k_1+\cdots+k_{j-1}+1}(\bM_N) - \rho_{\theta_j})=V_{1, N} +  o_{\mathbb P}(1)
$$
where 
$$V_{1,N} = U_{K_j \times 1}^* [B_{k,N}]^\nwarrow _{K_j}U_{K_j \times 1}. $$
It follows that $V_{1, N}$ converges towards the gaussian distribution ${\cal N}(0,\frac{t}{2}\frac{\theta_j^2 \sigma^2}{\theta_j^2- \sigma^2})$
and then according to Proposition \ref{convBinfinite}, 
  $\max_{i=1}^{K_j} \vert (U_k)_{i1}\vert $
converges to zero when $N$ goes to infinity. $\Box$\\

Let $\theta_j$ such that $\theta_j > \sigma$ and $k_j=1$. 
Let us  prove now the description given in subsection \ref{intermediate} of the fluctuations of $\lambda_{k_1+\cdots+k_{j-1}+1}(\bM_N)$ for some intermediate situations between Case a) and Case b).
 Let $m$ be a fixed integer number. Assume that for any $l =1,\ldots,m$
$(U_k)_{l 1}$ is independent of $N$,   whereas $\max_{m < l \leq K_j} \vert (U_k)_{l 1}\vert \rightarrow 0$ when $N$ goes to infinity. Following the proofs of Lemma \ref{premierbout} and Lemma \ref{deuxiemebout}, one can check that $V_{1, N}$ converges 
in distribution towards 
$\sum_{i,l=1}^m a_{il} \xi_{il}+ {\cal N}$ in the complex case, $\sum_{1\leq l \leq i \leq m} a_{il} \xi_{il}+ {\cal N}$
in the real case,
where $\xi_{il}, (i,l) \in \{1,\ldots,m\}^2, {\cal N}$ are independent random variables such that
\begin{itemize}
\item for any $(i,l) \in \{1,\ldots,m\}^2$, the distribution of $\xi_{il}$ is $\mu$;
\item $a_{il}=\left\{\begin{array}{lll}\sqrt{2} \Im (\overline{(U_k)_{l 1}}(U_k)_{i 1}) \mbox{~if~} i<l\\
\sqrt{t} \Re (\overline{(U_k)_{l 1}}(U_k)_{i 1}) \mbox{~if~} i>l\\ \sqrt{\frac{t}{2}}
\vert (U_k)_{l 1}\vert^2 \mbox{~if~} i=l; \end{array} \right.$
\item ${\cal N}$ is a centered  gaussian variable with variance 
$$ \frac{t}{4} \frac{\left[m_4 - 3 \sigma^4\right] \sum_{l=1}^m\vert (U_k)_{l 1}\vert^4}{\theta_j^2}   +
\frac{t}{2}\frac{\sigma ^4}{\theta_j^2-\sigma^2} + \frac{t}{2}\left[1- \left(\sum_{l=1}^m\vert (U_k)_{l 1}\vert^2\right)^2\right]\sigma^2.$$
\end{itemize}

Now, following the lines of Step 4 (using the results of Steps 1 and 2), we can conclude that 
$c_{\theta_j}\sqrt{N} (\lambda_{k_1+\cdots+k_{j-1}+1}(\bM_N) - \rho_{\theta_j})$ converges in distribution towards
the mixture of $\mu$-distributed or gaussian random variables $\sum_{i,l=1}^m a_{il} \xi_{il}+ {\cal N}$ in the complex case, $\sum_{1\leq l \leq i \leq m} a_{il} \xi_{il}+ {\cal N}$
in the real case.$\Box$

\section{Appendix}

In this section, we  recall some basic facts on matrices and some results  on random sesquilinear forms needed for the proofs of Theorems \ref{CLT} and \ref{CLTkN}.
\subsection{Linear algebra}
For Hermitian matrices, denoting by $\lambda_i$ the decreasing ordered eigenvalues, we have the Weyl's inequalities:
\begin{lemme}{(cf. Theorem 4.3.7 of \cite{HJ})} \label{Weyl}
Let B and C be two $N \times N$ Hermitian matrices. For any pair of integers
$j,k$ such that $1 \leq j,k\leq N$ and $j+k \leq N+1$, we have
$$\lambda_{j+k-1} (B+C) \leq \lambda_{j}(B) + \lambda_{k}(C).$$
For any pair of integers
$j,k$ such that $1 \leq j,k\leq N$ and $j+k \geq N+1$, we have
$$ \lambda_j(B) + \lambda_k(C) \leq \lambda_{j+k-N} (B+C).$$
\end{lemme}

In the computation of determinants, we shall use the following formula.
\begin{lemme}{\label{Lemma0}} (cf. Theorem 11.3 page 330 in \cite{BS4}) Let $A \in \mathcal M_{k}(\mathbb C)$ and $D$ be a
nonsingular matrix of order $N-k$. Let also $B$ and $~^tC$ be two
matrices of size $k\times (N-k)$. Then
\begin{equation} \label{calculdet}
 \det \left( \begin{array}{ll}
A ~ B \\ C ~ D
\end{array}\right)= \det(D) \det(A- B D^{-1} C).\end{equation}
\end{lemme}

\subsection{ CLT  for  random sesquilinear forms}
In the following, a complex random variable $x$ will be said {\it standardized} if $\mathbb{E}(x)=0$ and $\mathbb{E}(\vert x \vert^2)=1$.

\begin{theoreme}\label{BaiSilver98}(Lemma 2.7  \cite{BS1})
Let $B=(b_{ij})$ be a $N \times N$  Hermitian matrix and
$Y_N$ be a vector of size $N$ which contains i.i.d
standardized entries with bounded fourth
moment. Then there is a constant $K>0$ such that $$\mathbb E\vert  Y_N ^* B Y_N - {\rm{Tr}} B\vert^2 \leq K \mathbb \Tr (BB^*).$$
\end{theoreme}
\noindent This theorem is still valid if the i.i.d  standardized coordinates $Y(i)$  of $Y_N$ have a distribution depending on $N$ such that $\sup_N \E(|Y(i)|^4) < \infty$.

\begin{theoreme} \label{Silver} (cf. \cite{BY2} or Appendix by J. Baik and J. Silverstein in \cite{CDF} in the scalar case)
Let $A=(a_{ij})$ be a $N \times N$  Hermitian matrix and
$\{ (x_i, y_i), i \leq N \} $ a sequence of i.i.d centered vectors in $\C^K \times \C^K$ with finite  fourth
moment. We write $x_i = (x_{li}) \in \C^K$ and $X(l) = (x_{l1}, \ldots, x_{lN})^T$ for $1 \leq l \leq K$ and a similar definition for the vectors $\{Y(l), 1 \leq l \leq K\}$.  Set $\rho(l) = \E[\bar{x}_{l1}y_{l1}]$.  Assume that the following limits exist:
\begin{itemize}
\item[(i)] $\omega = \lim_{N \rightarrow \infty} \frac{1}{N}\sum_{i=1}^N a_{ii}^2$,
\item[(ii)] $\theta = \lim_{N \rightarrow \infty} \frac{1}{N} {\Tr} A^2 =  \lim_{N \rightarrow \infty} \frac{1}{N}\sum_{i,j =1}^N |a_{ij}|^2$,
\item[(iii)] $\tau  = \lim_{N \rightarrow \infty} \frac{1}{N} {\Tr} AA^T = \lim_{N \rightarrow \infty} \frac{1}{N}\sum_{i,j =1}^N a_{ij}^2$.
\end{itemize}
Then the $K$-dimensional random vector $\frac{1}{\sqrt N} \Big ( X(l) ^* A Y(l) - \rho(l) {\Tr}
A\Big ) $ converges in distribution to a Gaussian complex-valued vector  $G$ with mean zero. The
Laplace transform of $G$ is given by
$$\forall c \in \C^K, \quad \E[\exp(c^T G)] = \exp(\frac{1}{2} c^T Bc),$$
where the $K \times K$ matrix $B=(B(l,l'))$ is given by $B=B_1 + B_2 + B_3$ with:
\begin{eqnarray} \label{defBmatrix}
B_1(l,l') &=& \omega ( \E[\bar{x}_{l1}y_{l1}\bar{x}_{l'1}y_{l'1}] - \rho(l) \rho(l'))  \nonumber \\
B_2(l,l') &=& (\theta -\omega)\E[\bar{x}_{l1}y_{l'1}] \E[\bar{x}_{l'1}y_{l1}]  \\
B_3(l,l') &=& (\tau - \omega) \E[\bar{x}_{l1}\bar{x}_{l'1}] \E[y_{l1}y_{l'1}]. \nonumber
\end{eqnarray}
\end{theoreme}

\subsection{CLT for the empirical distribution of a Wigner matrix and applications}
\begin{theoreme}\label{BaiYao}(Theorem 1.1 in \cite {BY1}) Let
$f$ be an analytic function on an open set of the complex plane including $[-2
\sigma, 2\sigma]$. If the entries $((W_N)_{il})_{1 \leq i \leq l \leq N}$ of a general Wigner matrix $\bW_N$ of variance $\sigma^2$ satisfy the conditions
\begin{itemize}
\item[(i)] for $i \neq l$, ${\mathbb E}(\vert (W_N)_{il} \vert^4) = const$,
\item[(ii)] for any $\eta> 0$, $\lim_{ N \rightarrow + \infty} \frac{1}{\eta^4 n^2} \sum_{i,l}
{\mathbb E}\left[ \vert (W_N)_{il} \vert^4 \, 1\!{\sf I}_{\{\vert (W_N)_{il}\vert\geq \eta \sqrt{N}\}}\right]=0,$
\end{itemize}
then $ N \Big ( \tr_N(f(\frac{1}{\sqrt N} \bW_N))-\int f d \mu_{sc} \Big )$ converges in distribution towards a Gaussian variable, where $\mu_{sc}$ is the semicircle distribution of variance $\sigma^2$.
\end{theoreme}

\noindent We now prove some convergence results of the resolvent $\hat G$ used in the previous proofs.\\
Let $1 \leq j \leq J_{+ \sigma}$ and $k$ such that $\frac{k}{\sqrt{N}} \rightarrow 0$.
\begin{lemme}{\label{LemmaAppendix}} Each of the following convergence
holds in probability as $N \to \infty$:
\begin{itemize}
\item[i)] $ \sqrt{N} \left( \tr_{N-k} \hat{G}(\rho_{\theta_j}) - 1/\theta_j
  \right ) \longrightarrow 0$,
\item[ii)] $ \tr_{N-k} \hat{G}^2(\rho_{\theta_j}) \longrightarrow \int \frac{1}{(\rho_{\theta_j} - x)^2}d\mu_{sc}(x)= {1}/({\theta_j^2 - \sigma^2})$,
\item[iii)] $\frac{1}{N-k} \sum_{i=1}^{N-k} (\hat{G}(\rho_{\theta_j})_{ii})^2 \longrightarrow \left(\int \frac{d\mu_{sc}(x)}{\rho_{\theta_j} -x}\right)^2 =  {1}/{\theta_j^2}$.
\end{itemize}
\end{lemme}

\noindent {\bf{Proof of Lemma \ref{LemmaAppendix}: }} We denote by $G$ the resolvent of the non-Deformed
Wigner matrix ${W_{N-k}}/{\sqrt{N}}$. \\
\noindent i) By Theorem \ref{BaiYao}, one knows that $  \sqrt{N} \left( \tr_{N-k}
G(\rho_{\theta_j}) -  \int \frac{d\mu_{sc}(x)}{\rho_{\theta_j} -x} \right) $ converges in probability towards
0. Now, we have $\int \frac{d\mu_{sc}(x)}{\rho_{\theta_j} -x} = \frac{1}{\theta_j}$ (see \cite{HP} p. 94).  It is thus enough to show that
$$\tr_{N-k}\hat{G}(\rho_{\theta_j}) - \tr_{N-k}G(\rho_{\theta_j}) =o_{\mathbb P} ({1}/{\sqrt N}).$$
Let then $U_{N-k }:=U$ (resp. $D_{N-k}$) be a unitary (resp.
diagonal) matrix such that $ {A}_{N-k}= U^*D_{N-k}U$. Then, one has
\begin{eqnarray*} \vert \tr_{N-k}\hat{G}(\rho_{\theta_j})
- \tr_{N-k} G(\rho_{\theta_j}) \vert
 & = & \vert \tr_{N-k} \big (\hat{G}(\rho_{\theta_j}){A}_{N-k}G(\rho_{\theta_j}) \big
 )\vert\\
 & = & \vert \tr_{N-k} \big (D_{N-k} U^*
G(\rho_{\theta_j})\hat{G}(\rho_{\theta_j}) U \big ) \vert\\
& := & \vert \tr_{N-k} \big (D_{N-k} \Lambda (\rho_{\theta_j}) \big
)\vert \leq ({r}/{(N-k)}) \Vert D_{N-k} \Vert \Vert \Lambda
(\rho_{\theta_j}) \Vert
\end{eqnarray*}
where $r$ is the
finite rank of the perturbed matrix ${A}_{N-k}$. \\
One has
$\Vert D_{N-k} \Vert \leq \Vert A_{N} \Vert :=c$
(with $c=\max(\theta_1,|\theta _J|)$ independent from $N$). Moreover on the event $\tilde{\Omega} _N := \Omega^{(2)} _N \cap \{ \Vert
{W_{N-k}}/{\sqrt{N}} \Vert < 2 \sigma + \delta \}$, $\Vert \Lambda
(\rho_{\theta_j}) \Vert \leq (\rho _{\theta_j} - 2\sigma- \delta)^{-2} $
(use  \eqref{lem0}) so that we deduce that
$$ \vert \tr_{N-k} (\hat{G}(\rho_{\theta_j})) - \tr_{N-k}
(G(\rho_{\theta_j})) \vert 1\!{\sf I}_{\tilde{\Omega} _N} \leq
\frac{rc}{N-k} \, (\rho _{\theta_j} - 2\sigma -
\delta)^{-2}\to 0.$$
Using Theorem \ref{BaiYao} and the fact that
$\mathbb P
(\tilde{\Omega} _N )\to 1$, we obtain the announced result.\\

\noindent ii) It is sufficient to show that $\tr_{N-k}
\hat{G}^2(\rho_{\theta_j}) -\tr_{N-k} G^2(\rho_{\theta_j}) \to 0$ in
probability since, by Theorem \ref{BaiYao}, one knows that $
\tr_{N-k} G^2(\rho_\theta)$ converges in probability towards
$\int  \frac{1}{(\rho_{\theta_j} - x)^2}
d\mu_{sc}(x)$.\\
Using the fact that ${\rm{Tr}} (BC)={\rm{Tr}} (CB)$, it is not hard
to see that
\begin{eqnarray*} \tr_{N-k} \hat{G}^2(\rho_{\theta_j})
-\tr_{N-k} G^2(\rho_{\theta_j})
 & = & \tr_{N-k} \left ( ( \hat{G}(\rho_{\theta_j})+ G(\rho_{\theta_j}))
 (\hat{G}(\rho_{\theta_j})-G(\rho_{\theta_j})) \right )\\
 & = & \tr_{N-k} \left (\hat{G}(\rho_{\theta_j})   {A}_{N-k} G(\rho_{\theta_j}) (G(\rho_{\theta_j}) +\hat{G}(\rho_{\theta_j}) )  )
 \right )\\
 & = & \tr_{N-k} \left (D_{N-k} U
G(\rho_{\theta_j})  (G(\rho_{\theta_j})
+\hat{G}(\rho_{\theta_j}) )   \hat{G}(\rho_{\theta_j}) U ^* \right )\\
& := & \tr _{N-k} \big (D_{N-k} \Lambda' (\rho_{\theta_j}) \big )
\end{eqnarray*}
where the matrices $D_{N-k}$ and $U$ have been defined in i). We
then conclude in a similar way as before since on the event
$\tilde{\Omega} _N$, $ \Vert \Lambda' (\rho_{\theta_j})
\Vert \leq 2 (\rho _{\theta_j} - 2\sigma -\delta)^{-3 }$.\\

\noindent For point iii), we refer the reader \cite{CDF}. Indeed, it
was shown in Section 5.2 of \cite{CDF} that the announced
convergence holds in the case $k=1$ and for $G$ instead of
$\hat G$. It is easy to adapt the arguments of \cite{CDF} which
mainly follow from the fact that, for any $z \in \mathbb C$ such
that $\Im m (z)
>0$, $\frac{1}{N-k} \sum_{i=1}^{N-k} (\hat{G}(z)_{ii})^2$ converges
towards $g_\sigma^2(z) $. But this latter convergence was proved in
Section 4.1.4 of \cite{CDF}. $\Box$  \\

\noindent
{\bf Acknowledgments} We would like to thank the anonymous referees  for their  pertinent comments which led to an overall improvement of the paper.

\end{document}